\documentstyle{amlts}
\begin{document}
\input amssym.def
\input amssym.tex
\annalsline{153}{2001}
\received{March 16, 1998}
\revised{August 17, 2000}
\startingpage{533}
\def\bye{\AuthorRefNames [MGT]

 \end{document}}
 \font\tenrm=cmr10
%\input BoxedEps.Tex %for mac
%\SetTexturesEPSFSpecial %for mac
\input boxedeps.tex % unix
\SetepsfEPSFSpecial % unix
\HideDisplacementBoxes
\def\figin#1#2{\medbreak
$$
 {\BoxedEPSF{fig#1 scaled
#2}%
}%
$$
\medbreak\noindent}
\def\phigin#1#2{\medbreak
$$
 {\BoxedEPSF{plot#1 scaled
#2}%
}%
$$
\medbreak\noindent}

\def\ritem#1{\item[{\rm #1}]}
%--------------- Author macros ---------------
%for Bbb in amstex
\catcode`\@=11
\font\twelvemsb=msbm10 scaled 1100
\font\tenmsb=msbm10
%\font\ninemsb=msbm7 scaled 1100%msbm9
\font\ninemsb=msbm10 scaled 800
\newfam\msbfam
\textfont\msbfam=\twelvemsb  \scriptfont\msbfam=\ninemsb
  \scriptscriptfont\msbfam=\ninemsb
\def\msb@{\hexnumber@\msbfam}
\def\Bbb{\relax\ifmmode\let\next\Bbb@\else
 \def\next{\errmessage{Use \string\Bbb\space only in math
mode}}\fi\next}
\def\Bbb@#1{{\Bbb@@{#1}}}
\def\Bbb@@#1{\fam\msbfam#1}
\catcode`\@=12

 \catcode`\@=11
\font\twelveeuf=eufm10 scaled 1100
\font\teneuf=eufm10
\font\nineeuf=eufm7 scaled 1100%eufm9
\newfam\euffam
\textfont\euffam=\twelveeuf  \scriptfont\euffam=\teneuf
  \scriptscriptfont\euffam=\nineeuf
\def\euf@{\hexnumber@\euffam}
\def\frak{\relax\ifmmode\let\next\frak@\else
 \def\next{\errmessage{Use \string\frak\space only in math
mode}}\fi\next}
\def\frak@#1{{\frak@@{#1}}}
\def\frak@@#1{\fam\euffam#1}
\catcode`\@=12
%-------------- Author entries --------------------

\title{Ideal triangle groups, dented tori,\\ and numerical analysis}
\shorttitle{Ideal triangle groups} % Shortened version for headline title

  \acknowledgements{Supported by a
National Science Foundation research grant and an
Alfred P. Sloan research fellowship.} 
\author{Richard Evan Schwartz}

\institutions{University of Maryland, College Park, MD\\
{\eightpoint {\it E-mail address\/}: res@math.umd.edu}}

 \def\bp{|\underline{\overline{\times}}|}
\def\A{{\bf A}} 
\def\B{{\bf B}} 
\def\C{{\bf C}} 
\def\D{{\bf D}} 
\def\E{{\bf E}} 
\def\F{{\bf F}} 
\def\G{{\bf G}} 
\def\H{{\bf H}} 
\def\I{{\bf I}} 
\def\J{{\bf J}} 
\def\K{{\bf K}} 
\def\L{{\bf L}} 
\def\M{{\bf M}} 
\def\N{{\bf N}} 
\def\O{{\bf O}} 
\def\P{{\bf P}} 
\def\Q{{\bf Q}} 
\def\R{{\bf R}} 
\def\T{{\bf T}} 
\def\U{{\bf U}} 
\def\V{{\bf V}} 
\def\W{{\bf W}} 
\def\X{{\bf X}} 
\def\Y{{\bf Y}} 
\def\Z{{\bf Z}} 
\def\heartsuit{{\bf Z}} 

\phantom{cold}

\centerline{\bf Abstract} \bigbreak
We prove the Goldman-Parker Conjecture: 
A complex hyperbolic ideal triangle group is
discretely embedded in ${\rm PU}(2,1)$ if and only if
the product of its three standard generators
is not elliptic.  
We also prove that such a group is indiscrete
if the product
of its three standard generators is elliptic.
A novel feature of this paper is that it uses
a rigorous computer assisted proof to deal with
difficult geometric estimates.

\section{Introduction}

A basic problem in geometry and representation theory is
the {\it deformation problem\/}.
Suppose that $\rho_0: \Gamma \to G_1$ is a discrete embedding
of a finitely generated group $\Gamma$ into a Lie group $G_1$.
Suppose also that $G_1 \subset G_2$, where $G_2$ is a larger
Lie group. 
The deformation problem amounts to finding
and studying discrete embeddings 
$\rho_s: \Gamma \to G_2$, which extend $\rho_0$.

When $\Gamma$ is the fundamental
group of a surface, $G_1= {\rm Is} (\H^2)$, the isometry
group of the hyperbolic plane, and
$G_2={\rm Is}(\H^3)$, the isometry group of hyperbolic three space,
one is dealing with the
theory of quasifuchsian groups, which is quite well
developed.
(See, for instance,
the bibliography in
[T].)

The {\it complex hyperbolic plane\/} $\C\H^2$ is a $2$-complex-dimensional
manifold, which is negatively curved, K\"ahler, and also
a symmetric space.
It contains $\H^2$ as a totally real,
totally geodesic subspace, and is often considered
to be its complexification.
The
theory of deforming Is$(\H^2)$ representations into
Is$(\C\H^2)$ is also quite rich, though much less developed.
(For a representative sample of
such work, see
[FZ],
[GKL],
[GuP],
[KR], [Tol].)

\demo{{\rm 1.1.} The Goldman-Parker conjecture}
In [GP], Goldman and Parker took
one of the first steps on the road to 
a theory of complex hyperbolic quasifuchsian groups.
They defined and partially classified
the {\it complex hyperbolic ideal
triangle groups},  which are
representations 
$\rho_s: \Gamma \to {\rm Is}(\C\H^2)$.
Here $\Gamma$ is the free product
$\Z/2*\Z/2*\Z/2$.
The representation
$\rho_s$ maps the standard 
generators of $\Gamma$ to
distinct, order-two, {\it complex reflections}, 
such that any product of two distinct generators 
is parabolic.
A complex reflection is an isometry
which has, for a fixed point set, a
complex line in
$\C\H^2$.
(See \S 2.3 for more information.)

Modulo conjugation,
there is a one-parameter family
$\{\rho_s|\ s \in \R\}$ of such representations. 
The indexing parameter, $s$, is the tangent of
the {\it angular
invariant\/} of the ideal triangle formed by the
three complex lines fixed by the generators.
The angular invariant measures the extent to which
the vertices of the triangle fail to lie in a
totally real subspace.
(See \S 2.5 for a definition.)
The representation $\rho_0$ is the complexification
the familiar
{\it real ideal triangle group\/} generated
by reflections in the sides of an ideal
triangle in $\H^2$. 

Let $\underline s=(105/3)^{1/2}$ and 
$\overline s=(125/3)^{1/2}$.
According to Goldman-Parker:
 \enddemo
 
\nonumproclaim{Theorem {\rm [GP]}}
  If $|s|>\overline s$ then $\rho_s$ is not
a discrete embedding{\rm .}
If $|s| \leq \underline s$
then $\rho_s$ is a discrete embedding{\rm .}
\endproclaim

Here is a sketch of the first half of this result.
Let $g_s$ be the product of all three generators of
$\rho_s(\Gamma)$, taken in any order.
In [GP] it is shown that
$g_s$ is loxodromic for $|s| \in [0,\overline s)$,
parabolic for $|s|=\overline s$, and elliptic
for $|s|>\overline s$.
(See \S 2.1 for a classification of isometries.)
If $g_s$ is elliptic, with finite order,
then $\rho_s$ is not an embedding. If $g_s$ is
elliptic, with infinite order, then $\rho_s$
is not discrete. 

The theorem above is not sharp.
For $|s| > \underline s$ 
the analysis in [GP] breaks down,
but Goldman and Parker conjecture that $\rho_s$  remains
a discrete embedding 
for $|s| \in (\underline s,\overline s]$.
The significance of the result/conjecture combination is that it 
proposes the first complete description of a
complex hyperbolic deformation problem.

\demo{{\rm 1.2.} Results and methods}
The purpose of this paper is to prove a sharp version of
the Goldman-Parker conjecture. 
\enddemo

\nonumproclaim{Main Theorem}
 $\rho_s$ is a discrete embedding    if and only if
$g_s$ is not elliptic{\rm .}  Also{\rm ,} $\rho_s$ is
indiscrete if $g_s$ is elliptic{\rm .}
\endproclaim

Our indiscreteness proof uses some Galois
theory to show that
$g_s$ must have infinite
order when it is elliptic.

For our discreteness proof,
let $\Gamma_s=\rho_s(\Gamma)$.
Let $\partial \C\H^2$ be the ideal
boundary of $\C\H^2$.
We construct a 
surface-like set $\heartsuit(s) \subset \partial \C\H^2$, which we call a
{\it dented torus\/}. 
We then prove that
the orbit $\Gamma_s \heartsuit(s)$ 
consists of (essentially) disjoint
surface-like sets. 
This phenomenon feeds into a variant of
the Klein combination theorem to
prove that $\rho_s$ is a discrete embedding.

We found $\heartsuit(s)$ experimentally, using
our interactive program [S1].
The technical part of this paper consists
in making the pictures in [S1] into
rigorous proofs, via computation and
numerical analysis. 
While our proof is logically independent
from [S1], most of the ideas came
from analyzing the output of this program.

We use two kinds of
programs to verify our computational claims.
The first program, which computes
quantities depending on  the eigenvalues of $g_s$, 
is Mathematica code [M].
This code runs with much more precision than we need.
The second program, which is used for the
bulk of the computation, is C code [KeR].
Our C code uses only those
operations which conform to the IEEE
standards [I] and employs 
interval arithmetic to 
bound the roundoff errors rigorously.
Our code takes about 350 computational hours, when run
on a Sparc Ultra 5.   We ran the code, in parallel,
on 11 such machines, over the course of a weekend.

Our code includes a variety of interactive visual
displays, which allow the user to see, to a large extent,
that everything operates as claimed.
The code is available upon request.

We give explicit parametrizations
of the main objects in the paper,
as well as pseudocode expositions of the main
algorithms.  We hope this is enough information
for the inclined (and energetic) reader
to reproduce the results here
in an independent experiment.

\demo{{\rm 1.3.} Overview}
In Section~2 we review some basic complex hyperbolic geometry.

In Section~3 we define the triangle groups,
establish the indiscreteness part of our
main theorem, and
study the eigenvalues of $g_s$.

In Section~4, we study the action of $\Gamma_s$ on 
the {\it Clifford torus}, 
$$T=\{(z,w) \in \C^2|\ |z|=|w|=\sqrt 2/2\} \subset S^3.$$
(The ideal boundary of $\C\H^2$ can be naturally identified
with the three-sphere.)
We give an alternate proof of the Goldman-Parker discreteness result.
Though our proof breaks down
at the same place the proof in [GP]
breaks down,
it sheds new light on what should be done
in the critical interval $(\underline s,\overline s]$.
Our basic idea is to put some dents in the Clifford torus,
so that the resulting surface moves around better under
the action of the group.

In Sections~5 and 6 we construct the ``surfaces'' which replace the
Clifford torus.
In Section~5, we introduce a projectively natural coning operation,
the {\it hybrid cone construction\/}.
This procedure creates surfaces called
{\it hybrid sectors},   which are certain embedded solid triangles.
In Section~6 we introduce, for each parameter $s$, 
the dented torus, $\heartsuit(s)$, a surface
which is made from cutting out two solid triangles
from the Clifford torus
and gluing back  six 
hybrid sectors.  

In Sections~7 and 9 we give the bulk of the proofs.
In Section~7 we give the discreteness proof
of the Main Theorem for the parameter
$\overline s$, modulo two technical estimates,
 1 and 2.
In Section~8 we establish the Year Lemma, a result which shows
that we only need consider $365$ parameters in the
interval $J=[\underline{s},\overline{s}]$.
In Section~9 we give the general discreteness proof,
modulo the two estimates,
 1 and 2.

The remainder of the paper is devoted to establishing
and illustrating Estimates 1 and 2.
In Section~10 we describe the method by which
we estimate the location of a hybrid sector in space
using finitely many computations.
In Section~11 we detail the computation scheme
used to verify Estimates 1 and 2.
In Section~12 we discuss the implementation and computational
accuracy of our code.
We also give a record of the
actual calculations.
In Section~13 we illustrate
Estimates 1 and 2 with computer
plots, produced with our code.
\enddemo

\demo{Suggested reading}
We have indicated, through a series of remarks,
the material in Sections~2--6 which is strictly necessary for the
proof of the parabolic case.
Read only this material, and then
read Section~7.   Following this, skip to
Section~13, and look at the computer pictures, which
give compelling evidence
that Estimates 1 and 2 are true.
Afterwards, try to understand the general case.\enddemo

\demo{{\rm 1.4.} Related results}
This paper is the first in a series.
We prove in [S2] that the orbifold
at infinity, in the parabolic case,
is commensurable to the Whitehead link complement.
[S2] also gives a noncomputational
proof that there is some
$\varepsilon>0$ such that
$\rho_s$ is a discrete embedding for
$s \in (\overline s-\varepsilon,\overline s]$.

In [S3] we explain (among other things) how to build models 
out of string, of the limit sets of the complex
hyperbolic ideal triangle groups.

In [S4] we
consider an extremely deformed discrete representation
of the $(4,4,4)$-Hecke triangle group, and
use it to produce the first example of a
closed hyperbolic $3$-manifold which bounds
a complex hyperbolic $4$-manifold.

For some neat results and conjectures on the algebraic side of
the complex hyperbolic ideal triangle groups
see [Sa].   

For deformations of other triangle groups into
${\rm PU}(2,1)$, see [W-G].

For another example of an application of computation to
hyperbolic geometry, which is independent from ours,
see the monumental [GMT].
\enddemo

\demo{{\rm 1.5.} Acknowledgements}
I must say that writing this paper has been an
extremely gruelling experience, though sometimes an
exhilarating one.
I would like to thank Bill Goldman for his perspective and
encouragement as well as for the efficient proof
of Lemma \ref{flow}.
I would like to thank
Martin Bridgeman,  Peter Doyle,
Jeremy Kahn,  Josh Maher, Andy Mayer,
Robert Miner, 
Robert Meyerhoff, 
John Parker, Bill Thurston, and Justin Wyss-Gallifent for
interesting conversations relating to this work.
I would like to thank the tireless referees 
for many helpful stylistic suggestions.   Finally,
I would like to thank my wife Brienne Brown for her
C programming advice, and in particular for
suggesting the elegant subdivision algorithms
used in Section~10.8 and Section~11.1.

I dedicate this paper to my daughter, Lucina
Caroline Schwartz, who was born right around the time
I started thinking about triangle groups, and who is
now learning about triangles herself.
I also dedicate this
paper to the memory of Hanna Sandler, a great enthusiast for the
complex hyperbolic triangle groups, who passed away
this winter.
\enddemo

\section{Background}

2.1. {\it The complex hyperbolic plane}.
$\C^{2,1}$ is a copy of the vector space $\C^3$
equipped with the Hermitian
form
\begin{equation}
\langle u,v\rangle=u_1 \overline v_1+u_2 \overline v_2-
         u_3 \overline v_3. \end{equation}
The spaces $\C\H^2$ and $\partial \C\H^2$ are respectively the
projective images, in the complex projective plane $\C\P^2$, 
of
\begin{equation}
N_-=\{v \in \C^{2,1}|\ \langle v,v\rangle <0\}; \hskip 15 pt
N_0=\{v \in \C^{2,1}|\ \langle v,v\rangle=0\}
\end{equation}
(See [G,  p.67] or [E].) 
The map 
\begin{equation}
\Theta(v_1,v_2,v_3)=\bigg(\frac{v_1}{v_3},\frac{v_2}{v_3}\bigg)
\label{theta}
\end{equation}
takes $N_-$ and $N_0$ respectively to the open 
unit ball and unit sphere in $\C^2$.
Henceforth we identify 
$\C\H^2$ with the open unit ball.

Given a point $V \in \C\H^2 \cup S^3$, we will say that
$\tilde V \in \Theta^{-1}(V)$ is a {\it lift\/} of $V$ and
that a lift of the form
$(v_1,v_2,1)$ is {\it affinely normalized\/}.
We define the {\it vector\/} $\Theta^{-1}(V)$ as
the affinely normalized lift of $V$.

${\rm SU}(2,1)$ is the group of $\langle,\rangle$ preserving,
determinant $1$ complex linear transformations.
${\rm PU}(2,1)$ is the projectivization of ${\rm SU}(2,1)$, and
elements of ${\rm PU}(2,1)$ preserve $\C\H^2$.
Concretely, each $\tilde T \in {\rm SU}(2,1)$ determines an
element of ${\rm PU}(2,1)$ via
\begin{equation}
T=\Theta \circ \tilde T \circ \Theta^{-1}.
\end{equation}
The map $\tilde T \to T$ is a $3-1$ surjective
lie group homomorphism.

An element $T \in {\rm PU}(2,1)$ is called {\it loxodromic\/} if
$T$ has exactly two fixed points in
$S^3$, {\it parabolic\/} if it has exactly one fixed point
in $S^3$, and 
{\it elliptic\/} if it has a fixed point in
$\C\H^2$.   This classification is exhaustive and exclusive.
$T$ is called {\it ellipto-parabolic\/} if
$T$ is parabolic and also stabilizes a complex line in $\C\P^2$.
(See [G,  p.~203] for more details.)
For instance, the element $g_{\overline s}$
is ellipto-parabolic.  

$\C\H^2$ has a ${\rm PU}(2,1)$-invariant Riemannian metric, unique
up to scale.
Let $\delta$ denote the path metric induced by this Riemannian metric.
When normalized as in
[G,  p.~78],
\begin{equation} 
\cosh^2 \frac{\delta(v,w)}{2}=\frac{\langle \tilde v,\tilde w \rangle \langle
\tilde w,\tilde v\rangle}
{\langle\tilde v,\tilde v\rangle \langle\tilde w,\tilde w\rangle}
=\hat \delta(v,w).
\label{delta}
\end{equation}

\demo{{\rm 2.2.} Heisenberg space}
Let $\Re(z)$ and $\Im(z)$ denote the real and imaginary
parts of a complex number $z$.
The {\it Siegel domain\/} is the affine subset
\begin{equation}
\{(z,w)|\ \Re(w)> |z|^2\}  \subset \C^2 \subset \C\P^2.
\end{equation}
 
Let $\cal H$ stand for $\C \times \R$.
We call $\cal H$ {\it Heisenberg space\/}.
Given $p \in S^3$, 
a {\it Heisenberg stereographic projection from $p$\/}
is a map $\B: S^3-\{p\} \to \cal H$
which has the form $\B=\pi \circ \beta$.
Here $\pi(z,w)=(z,\Im(w))$ and $\beta$ is a complex projective
transformation of $\C\P^2$ which identifies
$\C\H^2$ with the Siegel domain.
We write $\infty=\B(p)$ in this case.
For instance
\begin{equation}
\B_0(z,w)=\bigg(\frac{w}{1+z},\Im \frac{1-z}{1+z}\bigg)
\label{cayley}
\end{equation}
is a Heisenberg stereographic projection which maps
$(-1,0)$ to $\infty$.
The map $\B_0$ is a close relative of the
{\it Cayley Transform\/} given in
[G,  p. 112].
The map $\B_0$ conjugates the ${\rm PU}(2,1)$ stabilizer
of $(-1,0)$ to certain affine maps of $\cal H$.

We call $E_0=\{0\} \times \R$ the
{\it center\/} of Heisenberg space.

The complex lines tangent to $S^3$ form a
canonical contact distribution on $S^3$.
We call this distribution $\cal E$.
Heisenberg stereographic projection
maps $\cal E$ to a corresponding
distribution in $\cal H$, which we give
the same name.
In $\cal H$, the distribution $\cal E$ is
the null distribution to the contact form
\begin{equation}
\omega=2ydx-2xdy+dt
\label{form}
\end{equation}
(see [G,  p.~124]).
In this formula, $\C$ has been identified with
$\R^2$ in the usual way.
$\cal E$ has
cylindrical symmetry.  It is preserved under
maps of the form
$$(z,t) \to (\lambda z,|\lambda|^2 t+ r); \hskip 15 pt
\lambda \in \C^*, \    r \in \R.$$
If $\Pi$ is a plane in $\cal E$, based at the point
$(z,t)$, then the maximum slope of a vector in
$\Pi$ is $2|z|$.

We say that a curve, in either $S^3$ or $\cal H$,
is $\cal E$-{\it integral\/} if its tangent
vector  at every point  is contained in $\cal E$.
The term $\cal E$-{\it integral\/} is sometimes
called CR {\it horizontal\/}.  We prefer to use
the former term, because the term
{\it horizontal\/} comes up elsewhere
in the paper.\enddemo

\vglue8pt

\demo{{\rm 2.3.} $\C$-circles}
A {\it complex slice\/} is the intersection of
a complex line in $\C^2$ with $\C\H^2$.   Complex slices are
totally geodesic subspaces, when considered as
Riemannian subspaces of $\C\H^2$.
A $\C$-{\it circle\/} (also known as a {\it chain\/}) is the intersection
of a complex line with $S^3$, provided this intersection is more
than a single point.
A $\C$-circle is a round circle,
transverse to $\cal E$.
A $\C$-{\it arc\/} is a nontrivial arc
of a $\C$-circle.  Given two points $p \not = q \in S^3$,
there is a unique $\C$-circle
containing $p$ and $q$.

Here is a method for parametrizing $\C$-arcs.
Suppose $X,Y,Z$ are three distinct points 
contained in a common $\C$-circle, $C$.
After switching the order of $X$ and $Y$ if necessary,
there are lifts $\tilde X,\tilde Y,\tilde X$,
unique up to multiplying all three vectors
by the same unit complex number, such that
\begin{equation}
\langle \tilde X,\tilde Y\rangle=i ; \hskip 15 pt
\langle \tilde X,\tilde Z\rangle=
\langle \tilde Y,\tilde Z\rangle=1
\end{equation}
(compare \S 2.5 below).
Note that $i\tilde X + t \tilde Z \in N_0$ for $t \in \R$, so that  
\begin{equation}
C_Z(X,Y;t)=\Theta(i \tilde X + t \tilde Z),\hskip 15 pt t \in [0,1],
\label{cpara1}
\end{equation}
parametrizes the $\C$-arc of $C$ which joins $X$ to $Y$ and which avoids $Z$.

Let $N_+=\C^{2,1}-N_--N_0$.
If $C$ is a $\C$-circle, then there is the {\it polar vector\/} 
$C^* \in N_+$, unique up to scaling, such that
$C=\{v \in N_0|\ \langle v,C^*\rangle=0\}$.
There is a unique involution
$I_C \in {\rm PU}(2,1)$ fixing $C$.
This map is computed by setting
$I_C=\Theta \circ I_{C^*} \circ \Theta^{-1}$, where
\begin{equation}
I_{C^*}(\tilde u)=-\tilde u+
\frac{2 \langle\tilde u,C^*\rangle}{\langle C^*,C^*\rangle}  C^*.
\label{reflection}
\end{equation}
(See [G,  p. 70].)
Another feature of the polar vector is this:
Two $\C$-circles, $C_1$ and $C_2$
are linked in $S^3$ if and only if
$\hat \delta(C_1^*,C_2^*)<1$, where
$\hat \delta$ is as in Equation~\ref{delta}.

We say that a {\it Heisenberg chain\/} is the image of a
chain under a Heisenberg stereographic projection.
Curves of the form $(\{z\} \times \R) \cup \infty$ are
chains. In particular, the center of $\cal H$ is a chain
(with $\infty$ deleted.)
All other Heisenberg chains are
ellipses which project to round circles under
the projection $\pi_{\C}: \cal H$$ \to \C$.
(See [G,  p.125].)
Let $C$ be such a chain, with center of mass $m$.
Let $E_m$ be the plane of $\cal E$ based at $m$.
Let $\hat E_m$ be the affine plane which is
tangent to $E_m$ at $m$.   We call $\hat E_m$
the {\it prolongation\/} of $E_m$.
We have $C \subset \hat E_m$.
We will say that a {\it round\/} Heisenberg
chain is one which is, itself, a round circle.
The center of mass of a round Heisenberg chain
is contained in the center of $\cal H$.
Such curves are
contained in planes of the form $\C \times \{t\}$.
\enddemo

\vglue8pt

\demo{{\rm 2.4.} $\R$-circles}
A {\it real slice\/}
is a totally real, totally geodesic subspace of
$\C\H^2$.    
Every real slice is isometric to the real slice
$\R^2 \cap \C\H^2$.
An $\R$-circle is the accumulation set,
in $S^3$, of a real slice.   Every $\R$-circle
is ${\rm PU}(2,1)$-equivalent to the particular
$\R$-circle $\R^2 \cap S^3$.
All $\R$-circles are $\cal E$-integral.

We say that an $\R$-circle is {\it round\/} if
it is a round circle.  
An $\R$-circle is round if and only if the real slice
containing it lies in a real $2$-plane in $\C^2$,
which is true if and only if the real slice contains $(0,0)$.
The geodesics in a round $\R$-circle are line
segments in $\C^2$.
Not all $\R$-circles are round.
An $\R$-{\it arc\/}
is a nontrivial arc of an $\R$-circle.
There is an $S^1$-family of $\R$-arcs containing
two points $p \not = q \in S^3$. 
The union of all such arcs is called a
{\it spinal sphere\/}.  (See [G,  p.152].)

Here is a method for parametrizing $\R$-arcs.
Suppose $X,Y,Z$ are three distinct points 
contained in a common $\R$-circle, $R$.
There are lifts $\tilde X,\tilde Y,\tilde X$,
unique up to scaling all three vectors
by the same complex number, such that
\begin{equation}
\langle \tilde X,\tilde Y\rangle=
\langle \tilde Y,\tilde Z\rangle=
\langle \tilde Z,\tilde X\rangle \in \R.
\end{equation}
(Compare \S 2.5 below.)
We call $\{\tilde X,\tilde Y,\tilde Z\}$ an {\it 
equalized choice\/} of lifts.
The curve
\begin{equation}
R_Z(X,Y;u)=\Theta((1-u) \tilde X + u \tilde Y+u(u-1) \tilde Z),
\hskip 15 pt u \in [0,1],
\label{rpara}
\label{rpara1}
\end{equation}
parametrizes the $\R$-arc of $R$ which joins $X$ to $Y$, but which
avoids $Z$.

We say that a {\it Heisenberg $\R$-circle\/} is the
image of an $\R$-circle under a Heisenberg stereographic
projection.
Any Heisenberg $\R$-circle $\gamma$ which
contains $\infty$ is (the extension of) a straight
line. 
We call these $\R$-circles {\it straight\/}. Note that
$\gamma$ has the form
$(L \times \{t\}) \cup \infty$ where $L \in \C$ is a line through the origin,
if and only if $\gamma$ is straight and
intersects the center of $\cal H$.
In this case we call $\gamma$ {\it level\/}. 
All other Heisenberg $\R$-circles are curves,
which project to lemniscates via the projection
$\pi_{\C}$.   (See [G,  p.139].)

Here is a geometric interpretation of our parametrization:
\proclaim{Lemma}
\label{aff}
Let $\B_0$ be as in Equation {\rm \ref{cayley}.}
If $Z=(-1,0)$ and $X,Y,Z$ lie on a common $\R$\/{\rm -}\/circle then
$\B_0(R_Z(X,Y;u))=(1-u) \B_0(X) + u \B_0(Y)${\rm .}
\endproclaim

\demo{Proof}
Since $\B_0$ maps the $\R$-circle containing $X,Y,X$ to
a straight line,
the point $P=\B_0(R_Z(X,Y;u))$ must be a convex combination
of $\B_0(X)$ and $\B_0(Y)$. 
Our parametrization of $\R$-arcs is natural in the following sense:
If, for $j=1,2$, the points
$X_j,Y_j,Z_j$ are contained on an $\R$-circle, then
an element of ${\rm PU}(2,1)$ which
carries $\{X_1,Y_1,Z_1\}$, in order,
to $\{X_2,Y_2,Z_2\}$, also carries 
$R_{Z_1}(X_1,Y_1;u)$ to $R_{Z_2}(X_2,Y_2;u)$.
Also, $\B_0$ conjugates the stabilizer of $Z=(-1,0)$
to affine maps of $\cal H$.
The facts just  mentioned imply that the particular
convex combination depends only on $u$, and not
on the triple of points chosen.

We will make the computation for the points
$X=(1,0)$, $Y=(0,1)$ and $Z=(-1,0)$; an
equalized choice  of lifts is
$\tilde X=(1,0,1)$,
$\tilde Y=(0,2,2)$ and
$\tilde Z=(-1,0,1)$.
Equation~13 gives
$$
R_Z(X,Y;u)=\bigg( \frac{1-u^2}{1+u^2},
\frac{2u}{1+u^2}\bigg).
$$
We have $\B_0(X)=(0,0)$ and
$\B_0(Y)=(1,0)$ and, from the calculation above,
$\B_0(P)=u$.   Since the result is true in
this case, it is true in all cases.
\enddemo

\demo{{\rm 2.5.} The angular invariant}
See [G,  pp.\ 210-214] for proofs of the results
in this section.

Given three points $X \not = Y \not = Z \in S^3$, the
{\it Cartan angular invariant\/} is 
\begin{equation}
\A(X,Y,Z)=
\arg [-\langle \tilde X,\tilde Y \rangle \langle \tilde Y,\tilde Z 
\rangle \langle \tilde Z,\tilde X \rangle]  \in [-\frac{\pi}{2},
\frac{\pi}{2}].
\end{equation}
We have $\A(X,Y,Z)=\pm \pi/2$ if and only if these points
lie in a common $\C$-circle, and 
$\A(X,Y,Z)=0$ if and only if these points
lie in a common $\R$-circle.\enddemo

\demo{{\rm 2.6.} The group law on $S^3$}
Recall that $S^3$ is a Lie group.
With $S^3$ as the unit sphere
in $\C^2$, as above, the group law on $S^3$ is
\begin{equation}
(z_1,w_1) \cdot (z_2,w_2)=
(z_1w_1-z_2 \overline w_2,
 z_1w_2 + z_2 \overline w_1).
\end{equation}
The right multiplication map
$R_X(Y)=Y \cdot X$ is an element of ${\rm PU}(2,1)$ which
acts isometrically on $S^3$, in the round metric,
and moves all
points by the same amount.  The point here is that
$R_X$ commutes with the complex structure, which is
left multiplication by $i$.\enddemo

\demo{{\rm 2.7.} The Clifford torus}
We do not know a reference for the material
in this section.

The {\it Clifford torus\/} is the torus
$T=\{(z,w)|\ |z|=|w|=\sqrt 2/2\} \subset S^3.$
Define $\arg: T \to T_0$ by
the formula
$\arg (z,w)=(\arg z,\arg w)$.
The preimages, under $\arg$, of curves of
slope $0$, $1$ and $\infty$ are all chains on $T$.
We call such\break chains, respectively,
{\it horizontal}, 
{\it vertical},  and {\it diagonal\/}.
The preimages,\break under $\arg$, of curves
of slope $-1$ are round $\R$-circles on $T$, having the form
$\{(z,w)|\ zw=\lambda\} \cap S^3$, where 
$|\lambda|=1.$

We say that a {\it Clifford triangle\/} is a
solid triangle in $T$ whose edges are
contained in chains, one from each of the three
foliations.   See Figure 3.1 for an example.

We say that a {\it Clifford rectangle\/} is a solid
rectangle whose boundary consists of two vertical
$\C$-arcs and two horizontal $\C$-arcs.
In parametrizing a Clifford rectangle we
wish to avoid the exponential function,
for computational reasons. 
Here is how we do it.
Let $Q$ be a Clifford rectangle.  Let
$H$ and $V$ be maximal horizontal and vertical $\C$-arcs in
$\partial Q$.

Let $H_1$ and $H_2$ be the endpoints of
$H$ and let $H_3$ be the midpoint
of $\check H-H$.   Here $\check H$ is
the $\C$-circle containing $H$.
Likewise define $V_1$, $V_2$ and $V_3$.
We parametrize $Q$ by
\begin{equation}
Q(t,u)=(\pi_1(C_{H_3}(H_1,H_2;t)),
        \pi_2(C_{V_3}(V_1,V_2;u))).
\label{rational}
\end{equation}
Here $\pi_1$ and $\pi_2$ are projections onto the
first and second factors, and the components
of our parametrization come from
Equation \ref{cpara1}.

$$ \BoxedEPSF{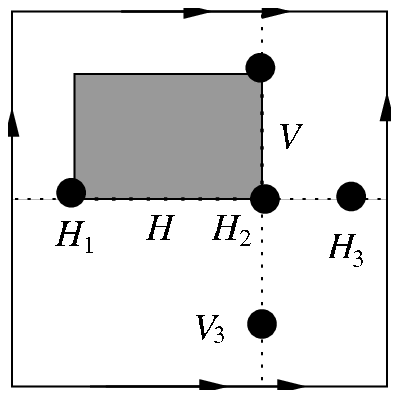 scaled 800}\qquad\quad\qquad$$
 \vglue-36pt
 \centerline{Figure 2.7}
  \vglue8pt

We call this the {\it rational parametrization\/}
of $Q$.
It appears that there are two choices for $H$ and two
for $V$.  However, recall from our definition of
$C_{H_3}(H_1,H_2;t)$ that the points
$H_1,H_2,H_3$ must be correctly ordered.
Likewise for $V_1,V_2,V_3$.
These extra constraints make the rational
parametrization unique.\enddemo

 2.8. {\it Pure translation flows}.

\demo{{R}emark}
The reader interested only in the parabolic case of the main
theorem may skip this section.
\enddemo

Say that a {\it pure translation\/} along a geodesic
$\gamma \in \C\H^2$ is one which preserves
all real slices containing $\gamma$.   Say that
a {\it pure translation flow\/} along $\gamma$ is
a flow $\phi_t: \C\H^2 \to \C\H^2$ such that
$\phi_t$ is a pure translation of $\gamma$, for
all $t$.
Let $I$ be the geodesic $(-1,1) \times \{0\}$.
Let $\|\ \|$ be the Euclidean norm on $\C^2$.

\proclaim{Lemma}
\label{flow}
The vector field $V(z,w)=(1-z^2,zw)$ generates a pure
translation
flow along $I${\rm .}
\endproclaim

\demo{Proof}
Let $\B_0$ be the map defined in Equation \ref{cayley}.
Let $(\zeta,\omega)=\B_0(z,w)$.
We compute
$$
d\B_0(V)=
\left[\matrix{-\frac{w}{(1+z)^2} & \frac{1}{1+z} \cr
\frac{-2}{(1+z)^2} & 0 } \right] \left[\matrix{z^2-1 \cr zw}\right]=
(\frac{w}{z+1},\frac{2-2z}{1+z})=(\zeta,2 \omega).
$$
This vector field generates the flow
$\psi_r(\zeta,\omega)=(r\zeta,r^2\omega)$.
Clearly, $\psi_r$ is a flow by complex linear maps which
preserves the
Siegel domain.  Hence, $\phi_t$ is a flow by
complex hyperbolic isometries, obviously preserving
$I$ and $\R^2$.
\enddemo

\proclaim{Lemma}
\label{flowbound}
Suppose that $\hat R$ is a real $2$\/{\rm -}\/plane which contains
$(0,0)${\rm .}  Suppose $p_1,p_2 \in \hat R-\C\H^2-S^3$ are
two points such that $(0,0),p_1,p_2$ are all contained on
a line $L \subset \hat R${\rm .}
Let $\beta$ be the pure translation
along $L \cap \C\H^2$ such that
$\beta(p_1)=p_2.$   Then{\rm ,} for all $q \in S_3${\rm ,}
$$\rho(q,\beta(q)) \leq \frac{\|p_1-p_2\|}{S^2-1};
\hskip 15 pt S=\inf(\|p_1\|,\|p_2\|).$$
\endproclaim

\vglue-20pt
{\it Proof}.
The stabilizer subgroup of $(0,0)$ in ${\rm PU}(2,1)$
acts isometrically on $S^3$ and transitively
permutes the real $2$-planes through $(0,0)$,
as well as the real lines through $(0,0)$.
Thus, we may assume that $\hat R=\R^2$,
that $L=I$, that
$p_1=(x_1,0)$ and $p_2=(x_2,0)$, with $1<x_1<x_2<\infty$.
Let $V$ be the vector field from Lemma \ref{flow}.
Since $\sup_{x \in S^3} \|V(x)\|=1$,  we conclude
that
$\rho(q,\phi_t(q)) \leq t$ for all $q \in S^3$.
The value of $t$, for which
$\phi_t(p_1)=p_2$, is
\medbreak
\hfill ${\displaystyle
t=\int_{x_1}^{x_2} \frac{dx}{x^2-1} \leq
\frac{x_2-x_1}{x_1^2-1}=\frac{\|p_1-p_2\|}{S^1-1}.
}$\hfill\qed

 \vglue12pt
\section{The triangle groups}
\advance \eqcount by 16

3.1. {\it Basic definitions}.
Given $s \in [0,\infty)$ let
$\Delta_s$ be the
Clifford triangle having vertices 
\begin{equation}
(\beta_s,\beta_s); \hskip 15 pt
  (\beta_s,\overline \beta_s); \hskip 15 pt
  (\overline \beta_s,\overline \beta_s); \hskip 30 pt
  \beta_s=\frac{(s+i)}{\sqrt{2+2s^2}}.
\end{equation}
The angular invariant of these three points
is $\pm \arctan s$, depending on the order in
which they are taken.
The first part of Figure 3.1 shows the images of these
points under the $\arg$ map.
 \figin{31.eps}{980}
 \vglue-36pt
  \centerline{Figure 3.1}
 \medbreak

Let $H_{+,s}$, $V_{+,s}$ and $D$ be the
horizontal, vertical, and diagonal $\C$-circles
containing pairs of these points.
Let $R$ be the $\R$-circle which
bisects $\Delta_s$.
Define
\begin{equation}
d(z,w)=(w,z); \hskip 15 pt
  r(z,w)=(\overline w,\overline z).
\label{dr}
\end{equation}
Note that $d$ is the complex reflection
fixing $D$ and $r$ is an antiholomorphic
involution fixing $R$.
Note that $d$ and $r$ commute.

We define $H_{-,s}=d(r(H_{+,s}))$ and
$V_{-,s}=d(r(V_{+,s}))$.
The second part of Figure 3.1 shows these chains.
Let $h_{\pm,s}$ be the complex reflection
fixing $H_{\pm,s}$, etc.
The group $\Gamma_s=\rho_s(\Gamma)$ is
generated by the triple
$\{d,h_{+,s},h_{-,s}\}$.  These are not the standard generators of $\Gamma_s$, but they are quite useful to us.
We define 
\begin{equation}
g_s=h_{+,s} \circ h_{-,s} \circ d.
\end{equation}
 Since $r$ conjugates $h_{\pm,s}$ to $v_{\pm,s}$ and
vice versa, 
\begin{equation}
r\Gamma_s r^{-1}=\Gamma_s; \hskip 15 pt
r g_s r^{-1}=g_s^{-1}.
\end{equation}
The representations $\rho_s$ and
$\rho_{-s}$ have the same image in
${\rm PU}(2,1)$, so we are justified in taking
$s \in [0,\infty)$.
We define $J=[\underline s,\overline s]$.
Usually we take $s \in J$.

\demo{{\rm 3.2.} Matrices}
We now give formulas for ${\rm SU}(2,1)$ lifts 
of some of the elements defined above.
\begin{equation}
\tilde h_{+,s}=\left[\matrix{-1&0&0 \cr
0&3&-4\overline \beta_s  \cr
0&4\beta_s  &-3} \right];
\hskip 15 pt
\tilde h_{-,s}=\left[\matrix{-1&0&0 \cr
0&3&-4\beta_s  \cr
0&4\overline\beta_s &-3} \right];
\end{equation}

\begin{equation}
\tilde d=\left[\matrix{0&-1&0 \cr -1&0&0 \cr 0&0&-1}\right].
\end{equation}

Setting $\tilde g=\tilde h_+ \tilde h_- \tilde d$,  we have
\begin{equation}
\tilde g_s=\left[ \matrix{0&-1&0 \cr
-A_1(s)&0&A_2(s) \cr
-A_2(s) &0&-\overline A_1(s)}\right]; \hskip 15 pt
 A_1(s)=\frac{s+17 i}{s+i}; \hskip 10 pt
A_2(s)=\frac{12 \sqrt 2 i}{\sqrt{1+s^2}}.
\label{As}
\end{equation}

Clearly, $\tilde d$ is a lift
of the map $d(z,w)=(w,z)$.
One can check easily that
$\tilde h_{\pm,s}$ is a determinant-one involution which
preserves $N_0 \subset \C^{2,1}$.  One can also check that 
this matrix has, as an eigenvector, a vector
which is polar to $H_{\pm,s}$.  One can also 
take traces of these matrices and their
products, checking that they agree 
with the traces of the corresponding
matrices in
[GP].\enddemo

\demo{{\rm 3.3.} Indiscreteness proof}
We now show that
$\rho_s$ is indiscrete for
$s \in (\overline s,\infty)$.
The work in [GP]
shows that
$g_s$ is elliptic for such  $s$.
We will show that
$g_s$ has infinite order.

From the formula for $\tilde g_s$, we see that
the trace $\tau_s$ of $\tilde g_s$ satisfies the
equation
\begin{equation}
\label{eq1}
(\tau_s+N+1)(\overline \tau_s+N+1)=N^2; \hskip 15 pt \tau_s \not = -1
\label{gal1}
\end{equation}
for $N=8$.   Note that
$\Re(\tau_s)<-1$ as long as $N \in \N$.
If $g_s$ has finite order then
\begin{equation} \label{eq2}
\tau_s=\omega_n^p+\omega_n^q+\omega_n^r;
\hskip 15 pt
\omega_n=\exp(2 \pi i/n); \hskip 15 pt p+q+r=0.
\label{gal2}
\end{equation}
Here $n$ is taken as small as possible.
Obviously $n>2$.
We suppose that Equations \ref{gal1} and \ref{gal2} are
true for some $N \in \N$, and derive a contradiction.
We set $\tau=\tau_s$.

\proclaim{Lemma}
\label{12}
$n$ divides $12${\rm .}
\endproclaim

{\it Proof}.
For $k$ relatively prime to $n$,
let $\sigma_k: \Q[\omega_n] \to \Q[\omega_n]$
be the Galois automorphism determined by
$\sigma_k(\omega_n)=\omega_n^k$.
Equation \ref{gal1} is defined
in $\Q[\omega_n]$, so it remains true if we replace $\tau$ by
$\sigma_k(\tau)$.   In particular,
$\Re(\sigma_k(\tau))<-1$.
Summing over all those $k \in \{1, \ldots ,n-1\}$ which are
relatively prime to $n$, we get
$\Re(\sum_k \sigma_k(\tau))<-\phi(n)$.
Here $\phi$ is the Euler phi function.
Hence \begin{equation} \label{eq3}
|\sum_k \sigma_k(\tau)|>\phi(n).
\label{gal3}
 \end{equation}
For $m=p,q,r$ let
$(m,n)$ be the greatest common divisor of $m$ and $n$.
Let $d_m=n/(m,n)$.  Note that  $\omega_n^m$ is a primitive
$d_m^{\rm th}$ root of unity, and the sum of all such roots is
one of $-1,0,1$.
The map $\Z/n \to \Z/d_m$ induces a map
$(\Z/n)^* \to (\Z/d_m)^*$, which is onto, with
multiplicity $\phi(n)/\phi(d_m)$.
Hence \begin{equation} \label{gal4}
|\sum_k \sigma_k(\omega_n^m)|
\leq \frac{\phi(n)}{\phi(d_m)}. \end{equation}
Combining Equations \ref{gal2}, \ref{gal3} and \ref{gal4} we get
\begin{equation} 
\label{gal5}
\frac{1}{\phi(d_p)} +
  \frac{1}{\phi(d_q)} +
  \frac{1}{\phi(d_r)} >1. \end{equation}
There is no positive integer $k$ such that $\phi(k)=3$,
and $\phi(k)=2$ only for $k \in \{2,3,4,6\}$.
We conclude that all three roots are $12^{\rm th}$ roots
of unity.
\hfill\qed

\demo{{R}emark}
For each choice of $N$, Lemma \ref{12} leaves only finitely
many cases to check.  One can finish the proof, for small
choices of $N$, by a quick computer search.
We did this for $N=8$, the case
of interest.   The argument following this remark finishes the
proof, in all cases, in a noncomputational manner.
\enddemo

\proclaim{Lemma}
\label{rp}
None of $p,q,r$ is relatively
prime to $n${\rm .} 
\endproclaim

\demo{Proof}
Suppose that $p$ is relatively prime to $n$.
All our congruences will be taken mod $n$.
Applying a Galois automorphism of
$\Q[\omega_n]$ we arrange that
$p=1$.   
Assume first that $n$ is even.
Since $\Re(\tau)<-1$ we must have
$\Re(\omega_n^q+\omega_n)<0$ and
$\Re(\omega_n^r+\omega_n)<0$.
This is only possible if
$q,r \equiv n/2$.   But then it is impossible to get
$p+q+r=0$.
Similarly, if
$n$ is odd then
$q,r \equiv (n \pm 1)/2$.
Since $p+q+r=0$, and $n>2$ we must have
$q \equiv r \equiv (n-1)/2$.
Here 
$\omega^q=\omega^r$ is a primitive $n^{\rm th}$ root of unity.
We apply a second Galois automorphism
to get $q=r=1$.
Thus $\Re(\omega^q+\omega^r)>0$,
a contradiction.
\enddemo

Lemma \ref{rp} says, in particular, that $n$
must have more than one prime factor.
Thus, $n \in \{6,12\}$.
Suppose $n=12$.
Since $\omega_n^p \omega_n^q \omega_n^r=1$,
we see that
$d_r$ divides the least common multiple $n'$ of
$d_p$ and $d_q$.   
Lemma \ref{rp} says $d_p,d_q,d_r<12$.
Suppose $n'<12$.  Since $d_p, d_q,d_r$ all
divide $n'$ we can replace $n$ by $n'$,
contradicting the minimality of $n$.
Thus, $n'=12$. 
Hence,
$(d_p,d_q)$ must be one of
$(3,4)$, $(4,3)$, $(4,6)$ or $(6,4)$.
From this it is easily seen that $d_r=12$,
a contradiction. 
The case $n=6$ is similar.

\demo{{R}emark}
The reader interested only in the parabolic case
can skip to the next chapter at this point.
\enddemo

\demo{{\rm 3.4.} Near monotonocity of the eigenvector}
Let $E_s$ be the $\C$-circle stabilized by the
element $g_s$, defined in Section~3.1.
Let $E_s^*$ be the affinely normalized vector which is
polar to $E_s$.
Note that $E^*_s$ is an eigenvector corresponding
to the norm-one eigenvector $\lambda_s$ of $g_s$.
Since the element $r$ conjugates $g_s$ to
$g_s^{-1}$, we have
$E^*_s=(e(s),\overline e(s),1)$.

Say that a function $f: J \to \R$ is $\varepsilon$-{\it monotone\/}
if there is a monotone function
$\tilde f: J \to \R$ such that $\|f-\tilde f\|_{\infty} \leq \varepsilon$.

\proclaimtitle{near monotonicity}
\proclaim{Lemma}
\label{near}
The functions $\arg e(s)$ and
$|e(s)|$ are $10^{-5}$\/{\rm -}\/monotone decreasing for $s \in J${\rm .}
\endproclaim

\demo{Proof}
We will show below that $|e'(s)|<3$ for all $s \in J$.
We explicitly evaluate the two functions above
for $120001$ maximally and evenly spaced points
in $J$, and observe that the two functions are
monotone decreasing on these points,
up to a roundoff error of $10^{-20}$. 
See Section~12 for details of the calculation.
Using the bound on the derivative, together with
the fact that $J$ has length less than $.54$, we get
the result of this lemma.

We turn now to our bound on $|e'(s)|$.
Let $A_1$ and $A_2$ be the functions defined in Section~3.2.
Let
$$X(s)=A_2(s) e(s) + \frac{\overline e(s)}{e(s)}.$$
We first establish the following implication:
\begin{equation}|e(s)|<1\ \&\ |X(s)|>\frac{4}{3} \hskip 15 pt
\Longrightarrow \hskip 15 pt |e'(s)|<3\ \&\  |X'(s)|<20.
\label{strap}
\end{equation}
Assume that $|X|>4/3$ and $|e|<1$ at some
parameter $s \in J$, which we suppress from
our notation.
Note that $$\tilde g(e,\overline e,1)=(\lambda_0 e,
\lambda_0 \overline e,\lambda_0)=
(- \overline e,-A_1e + A_2,-A_2e - \overline A_1).$$
Equating the first and third coefficients respectively gives
\begin{equation}
\lambda_0=-\frac{\overline e}{e}
\hskip 20 pt
e=\frac{-\overline A_1 - \lambda_0}{A_2}. \label{eq10} 
\label{theequation}
\end{equation}
Combining these equations we obtain
$A_2 e^2 + \overline A_1 e-\overline e=0$.
Differentiating with respect to $s$ we get
$$[A_2' e^2 + \overline A_1' e ]+
(2 A_2 e + \overline A_1) e' - \overline e'=0.$$
Basic calculus shows that
$|A_1'(s)|, |A_2'(s)|<.5$  for $s \in J$. 
If $|e|<1$ then the
bracketed expression above has norm less than $1$.
Therefore
\begin{equation}
|(2A_2 e + \overline A_1)e' - \overline e'|<1
\label{x} 
\end{equation}
Observe that
$$2A_2 e +\overline A_1=A_2 e + \frac{A_2 e^2 + \overline  A_1e}{e}=
A_2 e + \frac{\overline e}{e}=X.$$
Rearranging
Equation \ref{x} we arrive at
$$|e'(s)||X-u(s)| <1; \hskip 15 pt u(s)=e'(s)/\overline e'(s).$$
Since $u(s)$ is a unit complex number  and  
$|X|>4/3$ we see that $|X-u|>1/3$.  Hence  $|e'|<3$.
Note that $|A_2(s)|<3$ for $s \in J$.   Hence
$$|X'| \leq 2|A_2'||e|+2|A_2||e'|+|\overline A_1'| <1+18+.5<20.$$
We now know that Equation \ref{strap} is true. 

Say that $s \in J$ is a {\it good parameter\/} if
the two bounds on the left hand side of
Equation \ref{strap} hold at $s$.
Otherwise say that $s$ is a {\it bad parameter\/}.
Let $s_0, \ldots ,s_{55}$ be the set of $56$ maximally and evenly
spaced points in the interval~$J$.    
Note that $J$ has length $\sqrt{125/3}-\sqrt{35}<.55$.
Thus, every point of $J$ is within $.005$ of some point $s_j$.
As explicit calculation
shows that
$$
|X(s_j)| > 1.49; \hskip 15 pt 
|e(s_j)|<.9; \hskip 15 pt j=0,\ldots ,10.
$$
(See \S 12 for details of the calculation.)
Suppose $s \in J$ is some bad parameter.
Let $s_j$ be the point in our list nearest to $s$.
We can assume that all parameters in $[s_j,s]$ are
good.
This interval has length at most $.005$.
Equation \ref{strap},
says that 
$|X|$ and $|e|$ can each change by at most $.1$ on this interval.
Hence, $|X(s)| \geq 1.49-.1>1.3333...$ and
$|e(s)|<.8+.1<1$.
This is a contradiction.  Hence
every parameter is good.
Our bound follows immediately.
\enddemo
 
\vglue-6pt
\section{Images of the Clifford torus}
\advance\eqcount by 31
\vglue-6pt

\demo{{\rm 4.1.} A discreteness criterion}
We begin with a discreteness criterion.
This discreteness criterion is a close variant of
the ping-pong lemma and the Klein
combination theorem.  For original sources,
see [K] and [Mas].
A formulation which
is similar to the one here appears in
[Mac].

Recall that $\rho_s: \Gamma \to {\rm PU}(2,1)$ is our basic
representation.  We have also set $\Gamma_s=\rho_s(\Gamma)$.
Let $H_s \subset \Gamma_s$ be the group generated by
$h_{+,s}$ and $h_{-,s}$.
We say that a {\it compressing pair\/} for
$\rho_s$ is a pair of subsets
$(U_1,U_2)$, having the following properties:
First, 
$d(U_1)=U_2$.  Second, there is a proper subset
$V \subset U_1$ such that
$h(U_2) \subset V$ for all $h \in H_s$.

\proclaim{Lemma}  If $\rho_s$ has a compressing pair then
$\rho$ is a discrete embedding{\rm .}
\endproclaim

\vglue-12pt

{\it Proof}.
Any element $g\in\Gamma$ has the property that
$$\rho_s(g)= h_n\circ d\circ h_{n-1} \circ d\circ\dots\circ h_1$$
where $h_j\in H_s$. The elements $h_n$ and $h_1$ can be
trivial, but all the other $h_j$ are nontrivial. We perform induction
on $n$. The inductive hypothesis is the following: For $j=1,2$,
\medbreak
 1.  if $h_n$ is trivial, then $\rho(g)(U_j)\subset d(V)\subset U_2$;
\smallbreak 2.  if $h_n$ is nontrivial, then $\rho(g)(U_j)\subset h_n(U_2)\subset V$.
\medbreak\noindent 
When $n=1$, there is nothing to prove.

Suppose the induction hypothesis holds for $n-1$. If $h_n=1$, then
$$\rho(g)U_j  = \rho(d)\rho(h_{n-1}\circ d\dots d\circ h_1)  U_j
\subset d(V).$$
If $h_n\neq 1$, then
$$\rho(g)U_j =  \rho(h_n)\rho(d \circ h_{n-1}\circ d\dots d\circ h_1) U_j
\subset \rho(h_n)U_2 \subset V$$
 as desired.

As in the ping pong lemma, we have shown that
$\rho_s(g)$ does something uniformly nontrivial 
when $g$ is not the identity element.
Hence, 
$\rho_s$ is a
discrete embedding.
\hfill\qed

\demo{{\rm 4.2.} The projection}
Recall that $T$ is the Clifford torus.
$S^3-T$ consists of two open solid tori:
$$U_1=\{(z,w)|\ |z|<|w|\} \cap S^3; \hskip 15 pt
  U_2=\{(z,w)|\ |z|>|w|\} \cap S^3.$$
We will reprove the discreteness result in [GP]
by showing that $(U_1,U_2)$ is a compressing pair
for $\rho_s$ as long as $s \in [0,\underline s]$.
In this section we reduce this to a problem about
disks in the hyperbolic plane.

The condition $d(U_1)=U_2$ is true since
$d(z,w)=(w,z)$.
Let $W \subset \C\H^2$ be the closed disk
$\{z=0\} \cap (\C\H^2 \cup S^3)$. 
Let $\Pi_W: S^3 \to W$ be the projection
$\Pi_W(z,w)=(0,w)$.
For $j=1,2$, let
$W_j=\Pi_W(U_j)$.  By definition,
$W_2$ is the open disk of (Euclidean) radius
$\sqrt 2/2$ centered at $(0,0)$, and $W_1$
is the annular region $\{(0,w)|\ |w| \in (\sqrt 2/2,1]\}$.
\vglue-12pt
\figin{41.eps}{850}
\vglue-9pt
\centerline{Figure 4.2.}

\proclaim{Lemma}
\label{projj}
$(U_1,U_2)$ is a compressing pair
for $\rho_s$ if and only if 
$h(W_2)\break \subset W_1$
for all $h \in H_s${\rm .}
\endproclaim

{\it Proof}.
Using the formula for $h_+$ given in Section~2.3
we see that $h_{\pm} \circ \Pi_W=
\Pi_W \circ h_{\pm}$.  Hence,
\begin{equation}
h \circ \Pi_W=\Pi_W \circ h; \hskip 20 pt  
h \in H_s.
\end{equation}
From this we see that
$h(U_2) \subset U_1$ if and only if
$h(W_2) \subset W_1$.
If $h(W_2) \not \subset W_1$
then $(U_1,U_2)$ is certainly not a compressing pair.

Suppose $h(W_2) \subset W_1$ for all $h \in H$.
Note that $h(W_2)$ is a disk and
$W_1$ is an annulus.
Hence, every finite union of the form
$\bigcup_{h \in H'} h(W_2)$ is a proper
subset of $W_1$.
Here $H'$ is a finite subset of $H$.
The element $h_{\pm}|_W$ is an isometric rotation
about a point $p_{\pm} \in W$, as shown in Figure 4.2.
Hence, the product $h_+h_-$ is loxodromic, regardless of parameter.

The Euclidean diameter of the set
$h(U_1)$ tends to $0$ as the word length of
$h$ tends to $\infty$.   
This fact, together with the finite union
property, implies that
$X=\bigcup_{h \in H} h(W_2)$ is a
proper subset of $W_1$.   
Thus, $V=\Pi_W^{-1}(X)$ is a proper subset
of $U_1$
such that $h(U_2) \subset V$ for all $h \in H$.
\hfill\qed

\demo{{\rm 4.3.} Patterns of disks}
The disks $\overline W_2$ and $h_{\pm}(\overline W_2)$ are tangent at
the point $p_{\pm}$.   
Thus, the disks $hh_{\pm}(\overline W_2)$ and $h(\overline W_2)$ are tangent
at the point $h(p_{\pm})$.  The points of tangency
all lie on the geodesic $A$, through the points  $p_{\pm}$.  This geodesic is the translation axis for
$h_+h_-$.  
Figure 4.3 shows two pictures of the situation.
\vglue-24pt
\figin{ill42.eps}{225}  %fig42.eps
 \centerline{Figure 4.3.}
 \medbreak

\noindent 
In the second half of Figure 4.3, the complex slice
$W$ has been identified with the upper half plane
model and the geodesic $A$ has been identified with
the positive imaginary axis.

If $h$ has odd word length then
$h$ is an isometric rotation about a point not in $W_2$.
Hence, $h(W_2) \subset W_1$ for words of odd length.
If $h(W_2) \cap W_2=\emptyset$ for some
$h \in H$ having length $2k$ then $h'(W_2) \cap W_2=\emptyset$ for
all words $h'$ having longer (even) length.
The point is 
that the hyperbolic distance from the hyperbolic center of
$W_2$ to the hyperbolic center of $h'(W_2)$ exceeds
the hyperbolic distance from the hyperbolic center of
$W_2$ to the hyperbolic center of $h(W_2)$.
Our conclusion:  If $h(W_2) \cap W_2=\emptyset$, for
$h$ having word length $2$, then
$(U_1,U_2)$ is a compressing pair at
the parameter in question.  

By symmetry, $h_{+,s}h_{-,s}(W_2) \cap W_2=\emptyset$ if
and only if $h_{-,s}h_{+,s}(W_2) \cap W_2=\emptyset$.
A calculation, which we omit, shows that the 
distance between the hyperbolic center of
$W_2$ and the hyperbolic center of $h_{+,2}h_{-,s}(W_2)$
is monotone decreasing in $s$, and that the two disks are tangent
when $s=\underline s$.  
This establishes the
discreteness result in [GP].
\enddemo

\demo{{R}emark}
Our discreteness criterion parallels
[GP,  Lemma 3.2].
Our Lemma \ref{projj} parallels
[GP,  Lemma 3.3].
The analysis in this section parallels
the auxiliary lemma used to
prove [GP,  Lemma 3.3].
The omitted calculation parallels
the calculation in [GP]  
concluding their discreteness proof.
Our Clifford torus is contained
in the ideal boundary of the
Dirichlet polyhedron considered
in [GP].
The map $\Pi_W$ corresponds
to the map $\Pi_C$ used in
[GP].
\enddemo

 4.4. {\it Another point of view}.
The reader may think, at this point, that our proof above
fails by accident.  Let us give another point of view
which shows that the proof must fail for structural
reasons.   

The fixed point $p_{\overline s}$ of $g_{\overline s}$
is contained in the Clifford torus, by symmetry.
The point here is that the element $r$ preserves the
Clifford torus and conjugates $g_{\overline s}$ to its inverse.
Let $Q$ be the complex line tangent to $S^3$ at
$p_{\overline s}$.   The element $g_{\overline s}$
stabilizes $Q$, and rotates this plane nontrivially.
$Q$ is transverse to
the Clifford torus $T$.   
This means that $g_{\overline s}(T)$ must intersect $T$
transversely.   For $s$ sufficiently close to $\overline s$
the same phenomenon must occur, by continuity. 
One may think of
$\underline s$ as the cutoff for this phenomenon.
When $s \in (\underline s,\overline s]$,  the element
$g_s$ does not have the translational strength to compensate for
its rotational component in a direction transverse to $T$.
In short, $g_s$ does not have the strength to pull
$T$ off of itself.
 
\vglue-12pt
\section{The hybrid cone}
\advance\eqcount by 32

\demo{{\rm 5.1.} The basic construction}
Suppose that $(E,p)$ is a pair, where $E$ is a $\C$-circle
and $p \in E$ is a point.  

\proclaim{Lemma}
Suppose $x \in S^3-E${\rm .}   There is a unique $\R$\/{\rm -}\/circle
$\gamma=\gamma(E,p;x)$ such that $x \in \gamma$ and
$p \in \gamma$ and $\gamma \cap (E-p) \not = \emptyset${\rm .}
\endproclaim

\demo{Proof}
We normalize by a Heisenberg stereographic projection
so that $E=E_0$, the center of $\cal H$, and
$p=\infty$.
In this case, there is a unique level 
Heisenberg $\R$-circle containing $x$.
This is $\gamma(E_0,\infty;x)$.
\enddemo

\figin{51.eps}{900}
\centerline{Figure 5.1.}
\medbreak

Let $\Omega(E,p;x)$ be the arc
of $\gamma(E,p;x)$ which connects $X$ to $p$  and which is
disjoint from $E-p$. 
If $S \subset S^3-E$
we define
\begin{equation}
\Omega(E,p;S)=\bigcup_{x \in S} \Omega(E,p;x).
\end{equation}
We call $\Omega(E,p;S)$ the {\it hybrid cone\/}
of $S$ with respect to $(E,p)$.
We call the individual $\R$-arcs, in the union above,
the {\it foliating arcs\/}.

Our construction is natural:
$\Omega(T(E),T(p);T(S))=T(\Omega(E,p;S))$, when
$T$ is the restriction to $S^3$ of
an isometry of $\C\H^2$.

We define the {\it endpoint map\/}
\begin{equation}
f_{\Omega}(x)=\gamma(E,p;x) \cap (E-p).
\label{endpoint}
\end{equation}

\demo{{R}emark}
When $f_{\Omega}(S)$ is a single
point, 
$\Omega$ is a subset of a {\it spinal sphere\/}.
(See [G,  Ch. 5] for the theory of spinal spheres.)
Thus the hybrid cone construction generalizes
the spinal sphere construction. \enddemo

\demo{{\rm 5.2.} Hybrid disks and hybrid sectors}
Let $(E,p)$ be as above. 
If $S$ is a $\C$-arc, let $\widetilde S$ be the
$\C$-circle which contains $S$.
A {\it hybrid disk\/} is a
hybrid cone of the form $\Omega(E,p;\widetilde S)$.
A {\it hybrid sector\/} is a hybrid
cone of the form $\Omega(E,p;S)$, where $S$ is
a proper nontrivial $\C$-arc.
We will
sometimes write $\Omega(E,p;\widetilde S)=\widetilde \Omega(E,p;S)$.\enddemo

\proclaim{Lemma}
Let $\Omega=\Omega(E,p;S)$ be a hybrid sector.
$\Omega-\partial \Omega$ is an analytically embedded
open solid triangle{\rm .}  $\widetilde \Omega-p-C$ is
an analytically embedded punctured open disk{\rm .}
\endproclaim
\vglue-12pt
\figin{52.eps}{900}
\centerline{Figure 5.2.}

\demo{Proof}
We will consider the case of hybrid sectors,
the case of hybrid disks being similar.
We normalize  by a Heisenberg stereographic projection
so that $E=E_0$, the center of $\cal H$, and
$p=\infty$.
In this setting $\Omega(E_0,\infty;S)-\infty$ is a surface
with boundary, ruled by rays.
The ruling rays are arcs of level
Heisenberg
$\R$-circles.  The endpoints of the rays are contained
in $S$.  The linking condition implies that
the rays are pairwise disjoint.
The rays certainly vary analytically with points on $S$.
Our lemma is clear from this description. \phantom{forever}
\enddemo

When $\Omega(E,p;S)$ is normalized
as in the previous lemma, we say that
it is in {\it standard position\/}.
Figure 5.2 shows a picture of a standard
position hybrid sector.
The point $\infty$
is not shown.
If continued in the other direction, all the
foliating arcs of a standard position hybrid sector
intersect the center $E_0$.
This accounts for the radiating appearance
of the top view.

\demo{{\rm 5.3.} Canonical parametrization}
In this section we describe a parametrization for hybrid
sectors.  The formulas for the  parametrization are only important for
the proof in the sense that our code actually uses   them.

Let $(E,p;S)$ be as above.
Let  $q$ be the point of $E$ diametrically opposite~$p$.
(This makes sense because $E$ is a round circle.)
Let $\zeta \in S^3$ be some auxiliary
point chosen so that $p,q,\zeta$ lie on a common
$\R$-circle.   Let
$\{\hat p,\hat q,\hat \zeta\}$ be an equalized
choice of lifts.  
Let $E^*$ be the polar vector to $E_s$, chosen 
so that $\langle \zeta,E^*\rangle=
\langle \hat p,\hat q\rangle$.

Let $S$ be a $\C$-arc as above.
Let $S_1$ and $S_2$ be the endpoints of $S$,
and let $S_3$ be the midpoint of $\tilde S-S$.
Referring to Equation \ref{cpara1}
we define
\begin{equation}
S_t=C_{S_3}(S_1,S_2;t),
\hskip 15 pt
\hat S_t=\frac{\Theta^{-1}(S_t)}{\langle \Theta^{-1}(S_t),\hat p\rangle},
 \hskip 25 pt t \in [0,1].
\label{cpart}
\end{equation}
$\hat S_t$ is the lift of $S_t$ such that
$\langle \hat S_t,\hat p\rangle=1$.
Now define
\begin{equation}
\hat A_t=\Im(\langle \hat q,\hat S_t \rangle)\ \hat p+
          i  \hat q , \hskip 25 pt
A_t=\Theta(\hat A_t).
\end{equation}
Note that $A_t$ lies on the complex slice determined
by $p$ and $q$.  It is easily checked that
$\langle \hat A_t,\hat A_t \rangle=0$ and
$\A(p,S_t,A_t)=0$, so that $A_t=f_{\Omega}(S_t)$.
With reference to Equation \ref{rpara1},
the following map
is a parametrization of $\Omega(E,p;S)$:
\begin{equation}
\Omega(E,p;S;t,u)=R_{A_t}(p,S_t;u); \hskip 15 pt (t,u) \in [0,1]^2.
\end{equation}
\enddemo

 5.4. {\it The elevation map{\rm :}  Basics}.
With   notation as above, the formulas in this section
are only relevant to the proof in that we actually use
them in our code.

Recall that $\pi(z,w)=(z,\Im(w))$.   Define
\begin{equation}
Y_{\zeta}=\pi \circ \Theta \circ \hat Y_{\zeta}; \hskip 15 pt
\hat Y_{\zeta}(x)=   (\langle \tilde x,E^* \rangle,
              \langle \tilde x,\hat q \rangle, 
              \langle \tilde x,\hat p \rangle).
\label{elev}
\end{equation}

It is easy to see that the vector $V(z,w)=w\hat p + \hat q + w E^*$ is
null if and only if $\Re(w)=|z|^2$.    Evidently
$\Theta \circ \hat Y(V(z,w))=(z,w)$.
Thus $\Theta \circ \hat Y \circ \Theta^{-1}$ maps
$S^3-p$ into the boundary of the Siegel domain.
In short,
$Y_{\zeta}$ is a Heisenberg stereographic
projection mapping
$\Omega(E,p;S)$ to standard position, normalized so that
$Y_{\zeta}(p)=\infty$ and $Y_{\zeta}(q)=(0,0)$ and
$Y_{\zeta}(\zeta)=(1,0)$.
\smallbreak
Let $S^1$ be the unit circle in $\C$.
Let $\Upsilon=S^1 \times \R$.
We define 
\begin{equation}
\xi(z,t)=(z/|z|,t); \hskip 15 pt
\pi_{\R}(z,t)=t.
\end{equation}
Next, the
$\zeta$-{\it elevation map\/} 
$\Psi: S^3-\{p\} \to \Upsilon$ is defined as follows:
\[
\Psi(x)=\left\{ \begin{array}{ll}
     \xi \circ Y_{\zeta}(x) & \mbox{if $x  \in S^3-E$} \\
     \\
     S^1 \times \pi_{\R} \circ Y_{\zeta}(x)
      & \mbox{if $x  \in E-p$} .
    \end{array}
\right.
 \] 
Technically,
$\Psi$ is a map from $S^3-p$ to
subsets of $\Upsilon$.   This detail
is absorbed into the statements of our estimates.

\proclaim{Lemma}
\label{nat}
$\Psi(\Omega(E,p;S)-\{p\})=\Psi(S)${\rm .}
\endproclaim

\demo{Proof}
Let $x \in S$ be a point.
By construction, $Y(\Omega(E,p;x))$ is a
ray of a level $\R$-circle.  Thus,
$\xi$ is constant on $Y(\Omega(E,p;x))$.
\enddemo

Say that a {\it barrier\/} for $(E,p)$ is a
pair $(f,\Psi)$,  where $f: S^1 \to \R$ is a
continuous function.
Let
$\Lambda_f$ be the graph of $f$ in
$\Upsilon$.
Let $\Lambda_{1,f}$ and
$\Lambda_{2,f}$ be the two components of
$\Upsilon-\Lambda_f$.
We will take $\Lambda_{2,f}$ to be the component
above $\Lambda_f$, as shown
in Figure 5.3.   In this figure, $S^1 \times \R$ is
mapped into the plane via the map $(z,t) \to (\arg z,t)$.

\figin{54.eps}{1000}
\centerline{Figure 5.3.}
\medbreak

Lemma \ref{nat} immediately implies the following result: 

\proclaim{Lemma}
\label{disjo1}
Let $\Omega=\Omega(E,p;S)$ be a hybrid sector{\rm .}
Suppose that $\gamma \in {\rm PU}(2,1)$
stabilizes $(E,p)${\rm .}
Suppose $\Psi(\Omega) \subset \Lambda_{1,f}$ and
$\Psi(\gamma(\Omega)) \subset \Lambda_{2,f}${\rm .}
Then $\Omega \cap \gamma(\Omega)=\{p\}${\rm .}
\endproclaim

\demo{{R}emark}
At this point, the reader interested
in the parabolic case of the Main
Theorem can skip to the next chapter. \enddemo

5.5. {\it The elevation map\/{\rm :} Fine points}.

\proclaim{Lemma}
\label{monomap}
Suppose that $\alpha$ is an $\R$\/{\rm -}\/arc{\rm ,}
contained in an $\R$\/{\rm -}\/circle $\hat \alpha${\rm .}
If $\hat \alpha$ intersects $E$
in two points{\rm ,} neither of which
is $p${\rm ,} and
$\alpha$ intersects $E$ in one point{\rm ,} then
$\pi_{\R} \circ \Psi$ is monotone
on $\alpha${\rm .}
\figin{541.eps}{1000}
\centerline{\rm Figure 5.4.}
 \endproclaim

\demo{Proof}
We set $Y=Y_{\zeta}$, and recall that
$\Psi=\xi \circ Y$, as above. Also,
$\pi_{\C} \circ  Y(\hat \alpha)$ is a lemniscate
centered at the origin and
$Y(\hat \alpha)$ is an $\cal E$-integral lift
of this lemniscate.   The slope of this
lift is only $0$ at the two points which
cover $0$.
In particular, the slope of
$Y(\alpha)$ vanishes only at the endpoint.
(See [G,  p.141] for an explicit parametrization.)
Hence, the height function is monotone on
$Y(\alpha)$. \phantom{skinow}
\enddemo

\proclaim{{C}orollary}
\label{disjo2}
Let $\Omega(E,p;S)$ be a hybrid disk{\rm .}
Suppose $\gamma \in {\rm PU}(2,1)$ is such that
$\gamma(E)=E$ but $\gamma(p) \not = p${\rm .}
Suppose that 
$\Psi \circ \gamma(\partial \Omega) \subset \Lambda_2$ and
$\pi_{\R} \circ \Psi \circ \gamma(p)>\max f${\rm .}
Then $\Psi(\alpha) \subset
\Lambda_{2,f}$ for some foliating arc
$\alpha$ of $\gamma(\Omega)${\rm .}
\endproclaim

\demo{Proof}
The curve $\Psi \circ \gamma(\partial \Omega)$ is a simple
closed curve on $\Upsilon$ which, by hypothesis, is
contained in $\Lambda_{2,f}$.   This curve is a generator
for the homology of $\Upsilon$.  From  this it follows
that there is some $x \in \gamma(\partial \Omega)$ such that
$\pi_{\R} \circ \Psi(x)>
\max f$.
Let $\alpha$ be the $\R$-arc of $\gamma(\Omega)$ which connects
$x$ to $p$.   By Lemma \ref{monomap}, we have
$\pi_{\R} \circ \Psi(y)>\max f$ for all $y \in \alpha$.
This obviously implies that $\Psi(\alpha) \in \Lambda_{2,f}$. \phantom{skinow}
\enddemo

 5.6. {\it Tame hybrid disks}.
Say that a hybrid disk
$\Omega(E,p;\widetilde S)$ is {\it tame\/}
if the surface $\Omega-\partial \Omega-p$ is 
never tangent to the distribution $\cal E$.

\proclaim{Lemma}
\label{tamecrit}
Let $\Omega(E,p;\widetilde S)$ be a hybrid disk{\rm .}
Let $E^*$ and $S^*$ be polar vectors to
$E$ and $\widetilde S${\rm .}
If $\hat \delta(E^*,S^*)<1/4$ then
$\Omega(E,p;\widetilde S)$ is tame{\rm .}
\endproclaim

\demo{Proof}
Modulo ${\rm PU}(2,1)$-equivalence, there is a one-parameter family of inequivalent hybrid disks.
To get a complete family, we normalize so that
$$
E^*=(0,1,0); \hskip 15 pt S^*=(1/a,1,1); \hskip 15 pt
p=(-1,0).
$$
Here $a \in [0,1)$.
We compute $\hat \delta(E^*,S^*)=a^2$.
The map $\B_0$, defined in Equation \ref{cayley},
maps $\Omega$ to standard position.
The points
$$
p_1=(0,1); \hskip 15 pt 
p_2=\bigg(\frac{2a}{1+a^2},\frac{a^2-1}{1+a^2}\bigg)
$$
belong to $S$. 
By symmetry, $\B_0(S)$ intersects
$\R \times \{0\}$ at 
$\B_0(p_1)=(1,0)$ and
$\B_0(p_2)=(a-1,0)$ and has center of mass
$c=(a,0)$.
The point on $\pi_{\C} \circ \B_0(\Omega)$ having
minimum norm is $\pi_{\C} \circ \B_0(p_2)=a-1$.
If $a \in [0,1/2)$ then every point of
$\pi_{\C} \circ \B_0(\Omega)$ has norm greater
than $a$.

Let $T_x$ and $E_x$ be the tangent plane and $\cal E$-plane
at $x \in \B_0(\Omega-\widetilde S)-\infty$. 
Given any plane
$\Pi \in \cal H$,  we define the {\it slope\/}
$s(\Pi)= \sup_{v \in \Pi} \pi_{\R}(v)/|\pi_{\C}(v)|.$
From Equation \ref{form} we have
$s(E_x)=2|\pi_1(x)|>2|a|.$  Now,
$\B_0(\Omega)$ is ruled by horizontal rays, which
have endpoints on $\Sigma$, and which lie on lines
that intersect $\{0\} \times \R$.
This geometry implies that 
$s(T_x) \leq s(E_c)=2|a|<s(E_x)$.  Hence, $S_x \not = E_x$.
\enddemo

 \figin{53.eps}{1000}
 \centerline{Figure 5.7.}
 \medbreak

 5.7. {\it A disjointness criterion}.

\proclaim{Lemma}
\label{disjo3}
Suppose that $I \subset \R^n$ is some connected open set{\rm .}
Suppose{\rm ,}
for $t \in I${\rm ,} that
$\Omega(E,p_t;C_t)$ is a 
smoothly varying family of tame hybrid
disks and
$\Omega(E,q_t,S_t)$ is a smoothly
varying family of hybrid sectors or hybrid disks{\rm .}
%\smallbreak\noindent 
Suppose that
\medbreak
 {\rm 1.}  For each $t \in I${\rm ,}   $p_t \not = q_t${\rm .}
\smallbreak {\rm 2.}  $C_t \cap
\Omega(E,q_t;S_t)=\emptyset$ for all $t \in I${\rm .}
\smallbreak{\rm 3.}  $S_t \cap \Omega(E,p_t;C_t)=\emptyset$ for
all $t \in I${\rm .}
\smallbreak\hangindent=36pt\hangafter=1 {\rm 4.}  For some parameter $s \in I${\rm ,} one of
the foliating arcs of $\Omega(E,q_s,S_s)$
is disjoint from $\Omega(E,p_s,C_s)${\rm .}
\medbreak
\noindent 
Then $\Omega(E,p_t,C_t) \cap \Omega(E,q_t,S_t)=\emptyset$ for
all $t \in I${\rm .}
\endproclaim

\demo{Proof}
Let $N$ denote the union of all foliating arcs of
$\Omega(E,q_t,S_t)$.    The parameter $t$ ranges throughout $I$
in this definition.
We equip $N$ with the
Hausdorff topology:   Two elements of $N$ are
close if they are contained in small tubular
neighborhoods of each other.
Obviously, $N$ is connected.

Let $N' \subset N$ denote those arcs $\beta$
such that $\beta$ is a foliating arc for
$\Omega(E,q_t;S_t)$ and $\beta$ is disjoint from
$\Omega(E,p_t;C_t)$.   To establish this lemma,
it suffices to prove that $N'=N$.   Note that
$N'$ is nonempty, due to statement 4 above.
Clearly, $N'$ is  open.

To finish the proof we show that $N'$ is closed.
Let $\beta_1,\beta_2...$ be a sequence of elements
of $N'$ which converge to some $\beta \in N$.
Let $t$ be the parameter associated to
$\beta$.  
\medbreak
 1.  Statement 1
above implies that $\beta \cap p_t=\emptyset$.
\smallbreak
2.  Statement 2 implies that
$\beta \cap C_t=\emptyset$.
\smallbreak 3. Statement 3 implies that $\partial \beta \cap \Omega(E,p_t,C_t)=\emptyset$.
\medbreak\noindent 
The only possibility is that $\beta-\partial \beta$ intersects 
$\Omega(E,p_t,C_t)-C_t-p_t$.   If this intersection
is transverse, then nearby arcs, including some
of the $\beta_n$, would intersect nearby
hybrid disks. 
Hence $\beta-\partial \beta$ is tangent to
$\Omega(E,p_t;C_t)-C_t-p_t$ at some point $x$.

There is a 
foliating arc of $\Omega(E,p_t,C_t)$ which
contains $x$.  Call it $\alpha$.   Since $\Omega(E,p_t,C_t)$
is tame, the plane $E_x$ of $\cal E$
intersects the tangent plane $T_x$ to
$\Omega(E,p_t;C_t)$ transversely, in 
a line $L$. 
Both $\alpha$ and $\beta$ are tangent to $\Omega(E,p_t,C_t)$
at $x$, and both $\alpha$ and $\beta$ are $\cal E$-integral.
Hence, $\alpha$ and $\beta$ are both tangent to $L$.
The important conclusion here is that
$\alpha$ and $\beta$ are therefore tangent to each other.

We normalize so that
$\Omega(E,p_t;C_t)$ is in standard position.  Also,
$\pi_{\C}(\alpha)$ is a straight line segment,
contained on a ray through the origin and
$\beta$ is a Heisenberg $\R$-arc, contained in
a Heisenberg $\R$-circle $\hat \beta$.
The curves $\pi_{\C}(\alpha)$ and
$\pi_{\C}(\beta)$ are tangent at
$\pi_{\C}(x) \not = 0$.

There are two cases to consider.
Suppose first that
$\pi_{\C} (\beta)$ is a straight line.
In this case,
$\pi_{\C}(\hat \alpha)$ and $\pi_{\C}(\hat \beta)$
are the same line through the origin.
This implies
that $\hat \alpha=\hat \beta$.
Hence $\partial \alpha \cap \beta
\not = \emptyset$ or
$\partial \alpha \cap \beta \not = \emptyset$.
Either possibility contradicts our assumptions.

We conclude that $\pi_{\C}(\hat \beta)$ is a lemniscate.
Since $\hat \beta$ intersects $E_0$ twice,
$\pi_{\C}(\hat \beta)$ is   centered
at the origin.
However,
a straight line ($\pi_{\C}(\hat \alpha)$) through
the origin cannot be
tangent to a lemniscate centered at
the origin $(\pi_{\C}(\hat \beta)$) at a point ($\pi_{\C}(x)$)
other than the
origin.  This is a contradiction.
\enddemo

\section{Dented tori}
\advance\eqcount by 39

6.1. {\it Technical calculation}.

\demo{{R}emark}
For the reader only interested in the parabolic
case, the first two sections of this chapter can
be replaced by a single calculation at the
parameter $\overline s$.   One can just skim
this material. 
\enddemo 

Recall that $E^*_s=(e(s),\overline e(s),1)$ is the
eigenvector to $g_s$, as in Section~3.4.
Let $S$ be the circle of radius $\sqrt 2/2$.
Let $p_1(s)$ and $q_1(s)$ be the
points such that the tangent lines to $S$ at
$p_1(s)$ and $q_1(s)$ contain $e(s)$.
That is:
\begin{equation}
p_1(e)=\frac{1+ i \sqrt{2|e|^2-1}}{2 \overline e}, \hskip 20 pt
p_2(e)=\frac{1- i \sqrt{2|e|^2-1}}{2 \overline e}.
\end{equation}

Figure 6.1 shows a picture which is visually accurate
for all $s \in J$.

The following result implies that $p_1(s)$ never exits
the shaded region of the picture and $q_1(s)$ never
enters the shaded region.

\proclaim{Lemma}
\label{axiscalc}
For all $s \in J${\rm ,} $\arg p_1(s) \in [\pi/12,\pi]$ and
$\arg p_2(s)\in\break [\pi,2 \pi+\pi/12]${\rm .}
\endproclaim

\phantom{someday soon}
\vglue-1in
$$\BoxedEPSF{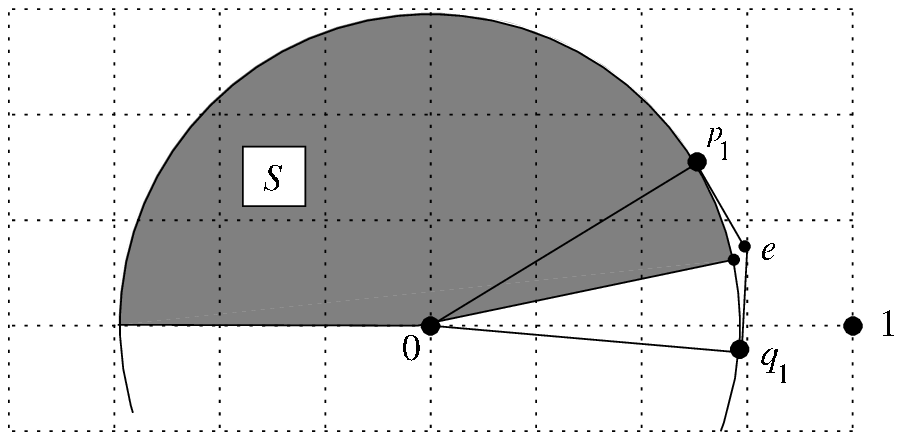 scaled 900}$$
\medbreak
 \centerline{Figure 6.1}

\demo{Proof}
Define
$\theta_1=\arg(e(\underline s))+10^{-5}$ and
$\theta_2=\arg(e(\overline s))-10^{-5}$.  Also define
$r_1=|e(\underline s)|+10^{-5}$ and
$r_2=|e(\overline s)|-10^{-5}$.  Finally, define
$e_{ij}=r_i \exp(\theta_j)$.
From the near monotonicity lemma,
and basic geometry,
$p_1$ lies in the interval bounded by the outermost points
$p_1(e_{ij})$ and
$q_1$ lies in the interval bounded by the outermost points
$q_1(e_{ij})$.
Explicitly calculating, we obtain this lemma.
Indeed, the points $p_1(e_{ij})$ and
$q_1(e_{ij})$ are positioned almost exactly as those in
Figure 6.1, independent of indices.
\enddemo

\demo{{\rm 6.2.} The axis lemma}
Define
\begin{equation}
\sigma=\tan(\pi/2-\pi/12)=2 + \sqrt 3.
\end{equation}
Let $\Delta_{\sigma}$ be the corresponding Clifford triangle,
on the Clifford torus $T$.
Let $\partial_h \Delta$ be the arc of the horizontal
$\C$-circle which bounds  $\Delta$.
Let $\partial_v \Delta$ be the arc of the vertical
$\C$-circle which bounds $\Delta$.
Let $\partial_d \Delta$ be the arc of the diagonal
$\C$-circle which bounds $\Delta$.   Thus, we
have
$\partial \Delta=\partial_h \Delta \cup
\partial_v \Delta \cup \partial_d \Delta$.
\enddemo

\proclaimtitle{axis lemma}
\proclaim{Lemma}
$E_s \cap T$ consists of two points{\rm .}
One of them{\rm ,} $p_s${\rm ,} is contained in
$\Delta_{\sigma}${\rm .}  The other point{\rm ,}
$q_s$ is not{\rm .}
Also $E_s$ links each of the
three $\C$-circles in $\partial \Delta_{\sigma}${\rm .}
\endproclaim

\demo{Proof}
We repeat the notation from
Section~6.1.
Let $S$ be the disk of radius $\sqrt 2/2$ centered
at the origin of $\C$.   Let $p_1(s)$ and
$q_1(s)$ be the two points of $\partial S$
such that the tangent lines through these
points contain $e(s)$.
Define
\begin{equation}
p_s=(p_1(s),\overline p_1(s)); \hskip 15 pt
  q_s=(q_1(s),\overline q_1(s)).
\end{equation}
We claim that
$E_s \cap T=p_s \cup q_s$.
By consruction, $p_s, q_s \in T$.
Recall that $R=\{z=\overline w\} \cap S^3$.
Note that $p_s, q_s \in R$.
Since $R$ is a round $\R$-circle, it
is contained in the real slice
$\hat R=\{z,\overline z\} \subset \C^2$.
The map $\hat R \to \C$ given
by $(z,\overline z) \to z$
carries $R$ to the disk $S$ and carries
geodesics in $\hat R$ to straight lines
in~$\C$.  Hence, the complex
lines tangent to $S^3$ at $p_s$
and $q_s$ both contain $E^*_s$.
This is enough to conclude that
$E^*_s$ is Hermitian perpendicular to
$p_s$ and $q_s$.
Hence, $E_s \cap R=p_s \cup q_s$.
If some other point
$o_s$ lies in $E_s \cap T$, then
so does $r(o_s)$ which is contained in the
same diagonal chain as  $o_s$.
This would force $E_s$ to be a diagonal
$\C$-circle, which it is not.
The result of
Section~6.1 is equivalent to the statement that
$p_s \in \Delta$
and $q_s \not \in \Delta$ 
for all $s \in J$.

Calculations show that
\begin{equation}
\hat \delta(E^*_{\overline s},\partial_d\Delta^*)=.4;
\hskip 15 pt
\hat \delta(E^*_{\overline s},\partial_h\Delta^*)=
\hat \delta(E^*_{\overline s},\partial_v\Delta^*)=
.042938\ldots \, .
\label{tameevidence}
\end{equation}
Here $\partial_d\Delta^*$ is a polar vector to
the $\C$-circle containing $\partial_d \Delta$.
Likewise for the other vectors.

Let $x$ be any one of
$d,v,h$.  The calculation above shows that
$E_{\overline s}$ links the
$\C$-circle containing $\partial_x\Delta$.
Since
$E_s$ never intersects the $\C$-circle containing $\partial_x \Delta$ for
$s \in J$, we see that
$E_s$ always links this $\C$-circle.
\enddemo

\demo{{\rm 6.3.} Main construction}
We define
\begin{equation}
\Omega_x(s)=\Omega(E_s,p_s;\partial_x \Delta); \hskip 20 pt
x=h,v,d. \end{equation}
For any subset $S \subset S^3$ let
$dS=d(S)$.  Here $d(z,w)=(w,z)$.  Next, define 

\begin{equation}
\Omega_{\Delta}(s)=
\bigcup_{x=h,v,d} \Omega_x(s);
\hskip 15 pt
\Omega(s)=\Omega_{\Delta}(s) \cup d\Omega_{\Delta}(s).
\end{equation}

By construction, $\Omega(s)$ is the union of $6$ hybrid sectors.
Next, we define
 \begin{equation}
\Xi=T-\Delta-d\Delta; \hskip 25 pt
\heartsuit(s)=\Omega(s) \cup \Xi.
\end{equation}
We call $\heartsuit(s)$ the {\it dented torus\/}.
$\heartsuit(s)$ has been constructed by cutting out
two Clifford triangles from $T$ and gluing
back in $6$ hybrid sectors.   
Figure 6.3 shows some of these sets schematically.

$$
\BoxedEPSF{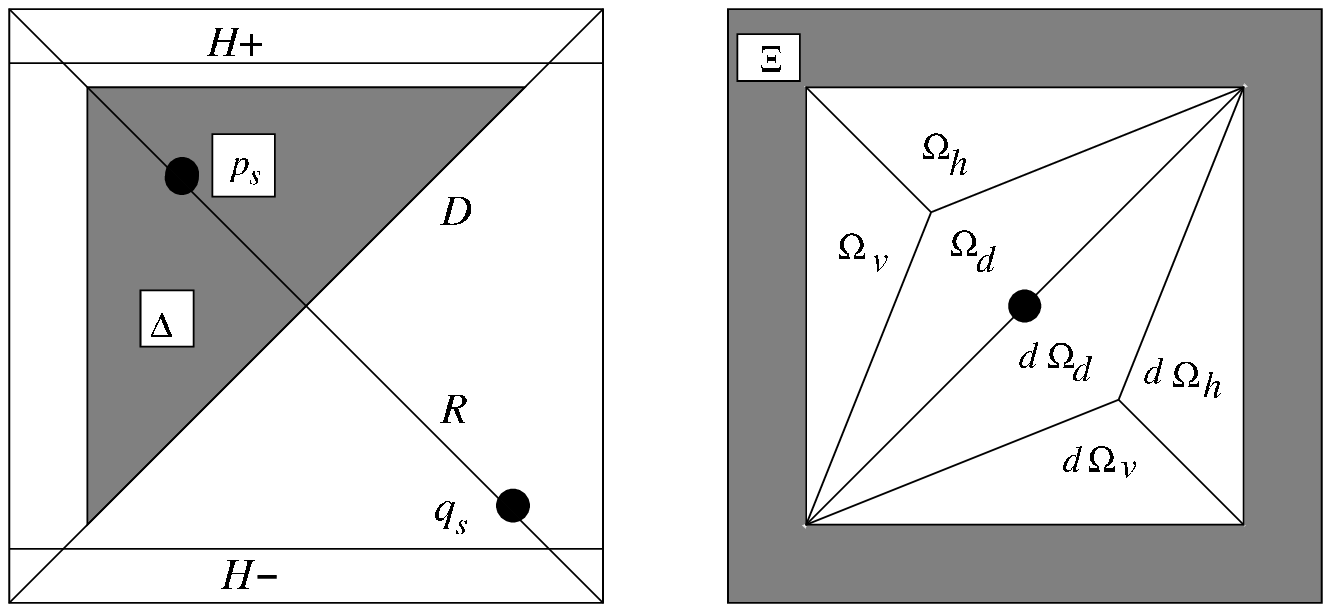 scaled 800}
$$
\vglue-12pt
 \centerline{Figure 6.3}
 \medbreak

$d$ preserves $\heartsuit(s)$ by construction.
$r$ preserves 
$\Delta$, $d\Delta$,
$\Xi$, $E_s$, $p_s$,
$dE_s$ and $dp_s$, all the sets which
are relevant to the construction.
Thus $\Z(s)$ has $4$-fold dihedral symmetry.
In particular, 
the points $p_s$ and $q_s$, which are
fixed by $r$, are diametrically
opposed on the circle $E_s$.  Hence, 
the canonical parametrization of
$\Omega_x(s)$, for $x=h,v,d$, uses the
point $q_s$ as the extra point.

Define
\begin{equation}
\zeta=\bigg(-\frac{1}{\sqrt 2},-\frac{1}{\sqrt 2}\bigg).
\label{zeta}
\end{equation}

Since $p_s,q_s, \zeta \in R$,
the $\zeta$-elevation map $\Psi_s$ is well defined for the
pair $(E_s,p_s)$.
By construction $\Psi_s(\Omega_x(s))$ is in standard position, for
$x=h,v,d$.
\enddemo

\section{The discreteness proof: Parabolic case}
\advance\eqcount by 47

7.1. {\it Five disks}.
Let $W$ be the complex slice used in Section~4.
Let $W^*$ be the complex line $\{z=0\}$ which contains
$W$ as its unit disk.
Let $W^*_+=\{(0,w) \in W^*|\ \Im(w) \geq 0\}$.
Define $W^*_-$ similarly.
Let $A_{\overline s} \subset W$ be the translation
axis of
the element $h_{+,\overline s}h_{-,\overline s}$.

We define five disks which play a key role in our proof.
\begin{itemize}
\item[1.] Let $W_A(\overline s) \subset W^*$ be the disk
such that $A_{\overline s} \subset \partial W_A(\overline s)$.
\item[2.] Let $W_B^+(\overline s)=h_{+,\overline s}(W_2)$.
\item[3.] Let $W_B^-(\overline s)=h_{-,\overline s}(W_2)$.
\item[4.] Let $W_C^+(\overline s)=h_{+,\overline s}(W_+*)$.
\item[5.] Let $W_C^-(\overline s)=h_{+,\overline s}(W_-*)$.
\end{itemize}

Let $W_{ABC}(\overline s)$ be the union of these five disks.
Figure 7.1 is a fairly accurate picture  (compare
Figures 13.1--13.6).  Henceforth we suppress the parameter $\overline s$.
 \figin{44.eps}{700}
\vglue-36pt
 \centerline{Figure 7.1}
  
Let $H_2$ be the set of words in $H$ which have length $2$.
\proclaimtitle{the five disk lemma}
\proclaim{Lemma}
\label{disjointcrit}
For all nontrivial $h \in H-H_2${\rm ,} 
$h(W-W_{ABC}) \subset W_{ABC}$ and
$h(W_2) \subset W_{ABC}${\rm .}
\endproclaim

{\it Proof}.
We first deal with $W-W_{ABC}$.
Suppose first $h$ has odd length.
By construction, $h$
preserves the axis $A$, and interchanges the
two components of $W-A$.
Therefore,
\begin{equation}
h(W-W_{ABC}) \subset h(W-W_A) \subset W_A
\subset W_{ABC}.
\end{equation}

By construction, $W-W_C^+-W_C^-$ is a fundamental
domain for the action of $h_+h_-h_+h_-$.
Thus, if $h$ has word length at least $4$ we have
\begin{equation}
h(W-W_{ABC}) \subset h(W-W_C^+-W_C^-) \subset
W_C^+ \cup W_C^- \subset W_{ABC}.
\label{39}
\end{equation}

We now turn to the images of $W_2$.
An explicit
calculation shows that $W_2 \cap W_C^{\pm}=\emptyset$.
Thus, $W_2 \subset W-W_C^{\pm} \subset W-W_{ABC}$.
Equation \ref{39} now implies that
$h(W_2) \subset W_{ABC}$ when $h$ has word length at least $4$.
If $h$ has word length $1$ then $h(W_2) \subset W_B$, by
definition.

Suppose $h$ has word length $3$.  By symmetry we can assume that
$h=h_+h_-h_+$.   Since $W_2 \cap W_B^{\pm}=\emptyset$
and $W_2 \cap W_C^{\pm}=\emptyset$, we have
$W_2-W_A \subset W-W_{ABC}$.   Thus, we are reduced to
showing that $h_+h_-h_+(W_2 \cap W_A) \subset W_{ABC}$.
Let $X=h_+(W_2 \cap W_A)$.  Let $X_{\pm}=X \cap W^*_{\pm}$.
By construction, $h_+h_-(X_+) \subset W_C^+$.
The region $X_-$ if nonempty, is contained in $W_A$,
by convexity, as shown in Figure 7.2.
Thus, 
$$h_+h_-(X_-) \subset h_+h_-(W_A)=W_A \subset W_{ABC}.$$
All in all, $h_+h_-h_+(W_2 \cap W_A)=h_+h_-(X) \subset W_{ABC}$, as desired.
\hfill\qed

\figin{45.eps}{1000}
\vglue-68pt
\centerline{Figure 7.2}
\eject

 7.2. {\it All words but two}.
Let $\Z=\Z(\overline s)$ be the dented torus
constructed in the previous chapter.
Recall that $\Z=\Omega \cup \Xi$.

\proclaimtitle{estimate 1: parabolic case}
\proclaim{Lemma}
$\Pi_W(\Omega_{\Delta}) \cap W_{ABC}=\emptyset${\rm .}
\endproclaim

{\it Proof}.
See Section~12.
\hfill\qed

\proclaim{{C}orollary}
\label{prom}
$\Pi_W(\Omega) \cap W_{ABC}=\emptyset${\rm .}
\endproclaim

{\it Proof}.
Let $r$ and $d$ be as in Equation \ref{dr}.
This corollary follows from the facts that
$\Omega=\Omega_{\Delta} \cup d\Omega_{\Delta}$,  
$rd$ commutes with $\Pi_W$,  
$rd(W_{ABC})=W_{ABC}$ and
$rd(\Omega_{\Delta})=d\Omega_{\Delta}$.
\hfill\qed

\proclaim{Lemma}
\label{pro}
If $h \in H-H_2$ is nontrivial then
$h(\Z) \cap \Z =h(\Xi) \cap \Xi.$
\endproclaim

{\it Proof}.
By the five disk lemma, and Corollary \ref{prom},
we have $h(\Pi_W(\Omega)) \subset W_{ABC}$.
Hence $\Pi_W(\Omega) \cap \Pi_W(h(\Omega))=\emptyset$.
Pulling back, we have $\Omega \cap h(\Omega)=\emptyset$.
By the five disk lemma, we have
$h(W_2) \subset W_{ABC}$.
Hence, by Estimate 1, we have
$\Pi_W(\Omega) \cap h(W_2)=\emptyset$.
Pulling back by $\Pi_W^{-1}$, and using the fact that
$\Pi_W(\Xi) \subset W_2$, we get
$\Omega \cap h(\Xi)=\emptyset$.
Since this is true for all nontrivial words $h \in H-H_2$
we have 
$\Omega \cap h^{-1}(\Xi)=\emptyset$.
That is,
$h(\Omega) \cap \Xi=\emptyset$.
This result now follows from the fact
that $\Z=\Omega \cup \Xi$.
\hfill\qed\medbreak

Let $x$ be the attracting fixed point of
$h_+h_-$. 
This point is one of the two intersection points of the axis
$A$ with $\partial W$.   In particular,
$$x=\Pi_W(x) \subset  W_A \subset W_{ABC}.$$
By Estimate 1, we have  $x \not \in \Z$.  Let $U_1$ be
the component of $S^3-\Z$ which contains $x$.

\proclaim{Lemma}
\label{struct2}
If $y \in W_{ABC}-W_2$ then
$\Pi_W^{-1}(y) \subset U_1${\rm .}
\endproclaim

{\it Proof}.
Note that $W_{ABC}-W_2$ is path connected.
Hence, there is a path $\delta_1 \subset W_{ABC}-W_2$
which joins $y$ to $x$, the fixed point of $h_+h_-$.
We can find a path $\delta_2 \subset S^3$
such that $\Pi_W(\delta_2)=\delta_1$.
One endpoint of $\delta_2$ is
necessarily $x$, since
$\Pi_W^{-1}(x)=\{x\}$.
We can make the other endpoint any point
of $\Pi_W^{-1}(y)$ we choose, since
the fibers of $\Pi_W$ are connected.
By construction,
$\delta_2 \subset S^3-\Z$.
\phantom{snow}\hfill\qed

\proclaim{{C}orollary}
\label{extra}
Both fixed points of $h_+h_-$ are contained in $U_1${\rm .}
If $\sigma$ is a fixed point of
$dh_+h_-d$ and $h \in H$ then
$h(\sigma) \in U_1${\rm .}
\endproclaim

{\it Proof}.
For the first claim, the fixed points
of $h_+h_-$ project to $\partial W \cap W_A
\subset W_{ABC}-W_2$.
For the second claim,
the fixed points of $h_+h_-$ have
zero first coordinate.  Hence, $\sigma$ has
zero second coordinate:
$\Pi_W(\sigma)=(0,0) \not \in W_A$.
If $h$ has odd length then
$\Pi_W(h(\sigma))=h(0,0) \subset W_A-W_2$.
If $h$ has even length, at least $4$, then
$\Pi_W(h(\sigma))=h(0,0) \subset W_C^{\pm}$,
which is disjoint from $W_2$.
If $h$ has length $2$
then $h(0,0)$ is contained in
$\partial W_C^{\pm}$.
We have already seen that
$W_2$ is disjoint from this
disk.  So, in all cases,
$h(0,0) \in W_{ABC}-W_2$.
Pulling back by $\Pi_W$, and applying
Lemma \ref{struct2} finish  the proof.
\hfill\qed

\proclaim{Lemma}
\label{distinct}
$U_1 \not = U_2${\rm .}
\endproclaim

{\it Proof}.
Let $y$ be the fixed point of the parabolic element
$h_+d$.
Recall that $H_+$ is the horizontal chain on the
Clifford torus which contains the horizontal
arc bounding $\Delta$, the Clifford triangle used
in the construction of $\Z$.
By construction $y \in H_+$.  Hence,
$\Pi_W(y)$ is the fixed point $p_+ \in W_2 \cap W_A$ of $h_+$.
In particular,
$\Pi_W(y) \subset W_{ABC}$.

Let $x$ be the fixed point of $h_+h_-$.
There is a path $\delta_0$, connecting $x$ to $y$, such that
$\delta_0-\{y\} \in U_1$, and $\delta_0$ is transverse to
$\Xi$ at $y$.  
The path $\delta_1=\delta_0 \cup d(\delta_0)$
joins $x \in U_1$ to $d(x) \in U_2$.
By construction, $\delta_1$
intersects $\Xi$ transversely in a single
point and is disjoint from
$\heartsuit-\Xi=\Omega$.
If $U_1=U_2$ then there is some path
$\delta_1'$ which connects $x$ to
$d(x)$ but which 
avoids $\heartsuit$.  The path
$\gamma=\delta_1 \cup \delta'_1$ is a closed
loop in $S^3$ which intersects $\Xi$ once, transversely,
and is disjoint from $\Omega$.
Hence $\gamma$ and
$\partial \Omega=\partial \Xi$ are algebraically linked.

Homological considerations say that
$\gamma$ links $\partial \Delta$ or
$\gamma$ links $d\partial \Delta$, or both.
Suppose, without loss of generality, that
$\gamma$ links $\partial \Delta$.
Since $\Omega_{\Delta}$ is the continuous
image of a disk, attached to $\partial \Delta$
by a degree one map, we see that
$\gamma$ and $\Omega_{\Delta}$ have nontrivial
intersection.  This is a contradiction. 
\hfill\qed 

\proclaim{Lemma}
\label{top}
Let $h \in H-H_2$ be nontrivial.
The component $U_1$
contains all
but one component of
$S^3-h(\heartsuit)${\rm .}
\endproclaim

{\it Proof}.
By Lemmas \ref{prom} and
  \ref{pro}, the intersection
$X=h(\heartsuit) \cap \heartsuit$ is
either the empty set, or one of
the two $\C$-circles $H_{\pm}$.
In either case,
$X$ neither disconnects $\heartsuit$
nor
$h(\heartsuit)$.
Basic topology now implies that
one component of $S^3-\Z$, which we call the
special component, contains all
but one component of $S^3-h(\Z)$ (and
{\it vice versa\/}).

We now show that
$U_1$ has nontrivial intersection
with more than one component of
$S^3-h(\heartsuit)$, meaning
that $U_1$ is the special component.
Since $h$ interchanges the fixed points of
$h_+h_-$, both of which are contained in
$U_1$, we see that $h(U_1) \cap U_1 \not =\emptyset$.
From Corollary \ref{extra}
we see that
$U_1$ also contains
$h(\sigma)$. 
Thus $h(U_2) \cap U_2 \not =\emptyset$.
\hfill\qed

\proclaim{Lemma}
 For all nontrivial
$h \in H-H_2${\rm ,} $h(U_2) \subset U_1${\rm .}
\endproclaim

{\it Proof}.
We already know that $U_1$ contains all
components of $S^3-h(\heartsuit)$, except for one.
We just have to rule out the possibility
that $h(U_2)$ is this one exceptional
component of $S^3-h(\Z)$.  Assuming to the contrary,
we get $h(U_1) \subset U_1$.
The element $rd$ conjugates $h$ to
$h^{-1}$ and preserves $\heartsuit$.
Furthermore, since $rd$ interchanges
the fixed points of $h_+h_-$, we
see that $rd(U_1)=U_1$.   
If $h(U_1) \subset U_1$ then
$$h^{-1}(U_1)=(rd)h(rd)(U_1)=
(rd)h(U_1) \subset rd(U_1)=U_1.$$
This implies that $h^{-1}(U_1) \cap U_2=\emptyset$.
By Corollary \ref{extra} we have
$h^{-1}(U_1) \cap U_2 \supset \{\sigma\} \not = \emptyset$,
a contradiction.
\hfill\qed
\medbreak

{\it Remark}.
Computer evidence suggests that
$\Z(\overline s)$ is a tamely embedded
topological torus and, in particular, that
$S^3-\Z(\overline s)$ consists of exactly
two connected components.  If this is true,
our argument above could be somewhat
simplified.  While it is
possible to convert the computer evidence into
a proof,  it is quite difficult to do so.
We settle for the slightly complicated arguments
above to avoid the extremely complicated arguments
needed to analyze the topology of $\Z(\overline s)$.

\demo{{\rm 7.3.} Words of length $2$}
Recall that $\Upsilon=S^1 \times \R^2$ and
$\Psi:S^3-p \to \Upsilon$ is the 
$\zeta=(1/\sqrt 2,1/\sqrt 2)$
elevation map associated to $(E,p)$.
Define the barrier function:
\begin{equation}
\label{specialf}
f(u)=t+.21\Im(\overline z u);
\hskip 20 pt
(z^2,2t)=\Psi \circ g(\zeta).
\end{equation}
In defining $z$, the
branch of $\sqrt{}$ is chosen so that
$0<\arg(z)<\arg(z^2)$.
Let $\Lambda$ and $\Lambda_j$ be sets associated
to the barrier $(f,\Psi)$.
We found this function experimentally.

\proclaimtitle{estimate 2: parabolic case}
\proclaim{Lemma}
For the parameter $\overline s${\rm :}
\smallbreak
{\rm 1.} $\Psi(\Xi \cup d\Omega_{\Delta}) \subset \Lambda_1(f)${\rm .}
\smallbreak {\rm 2.} $\Psi \circ g(d\Omega_{\Delta} \cup \Xi) \subset \Lambda_2(f)${\rm .}
\endproclaim

{\it Proof}.
See Section~12.
\hfill\qed

\proclaim{{C}orollary}
For $h \in H_2${\rm ,} $h(U_1) \subset U_1.$ 
\endproclaim

{\it Proof}.
It suffices to consider the element $h_+h_-$, by symmetry.
Note that $h_+h_-(\Z)=h_+h_-d(\Z)=g(\Z)$.
From Estimate 2 and Lemma \ref{disjo1} we have
$g(\Z) \cap \Z=\{p\}$.
Hence $h_+h_-(\heartsuit) \cap \heartsuit = \{p\}$.
Since $p$ disconnects neither
$\heartsuit$ nor $h_+h_-(\heartsuit)$, the same
argument as above shows that
$h_+h_-(U_2) \subset U_1$.
\hfill\qed\medbreak

All in all $(U_1,U_2)$ is a compressing pair
for $\rho$, and $\rho$ is a discrete embedding.

\vglue-6pt
\section{The year lemma}
\advance\eqcount by 50
\vglue-6pt

Recall that $J=[\underline{s},\overline{s}]$.  We  partition   $J$ into $365$
intervals of equal size.  The
set $J_0$ of midpoints of these intervals 
is   given by
\begin{equation}
(1-\frac{2k-1}{730}) \underline s +
\frac{2k-1}{730} \overline s; \hskip 20 pt k=1, \ldots ,365.
\end{equation}
Let $I_k$ be the partition interval containing $s_k$.

The geometry of $\Z(s)$ varies slowly with $s \in J$.  However,
this variation is extremely difficult to estimate directly.
In this chapter we replace the family of sets $\{\Z(s)|\ s \in J\}$
by a family of sets
$\{\Z(s_k,t)|\ t \in I_k; \hskip 5 pt k=1, \ldots ,365\}$.
The set $\Z(s_k,t)$ is a replacement for
$\Z(t)$, whose geometry is essentially the same as the
geometry of $\Z(s_k)$.
In this way, we only have to deal with the geometry of
the surfaces $\Z(s_k)$, for $k=1, \ldots ,365$.

\demo{{\rm 8.1.} Compatibility maps}
Let $s=s_k$, and $t \in I_k$ be some other parameter.
Let $(E_s,p_s)$ be the pair used to define $\Z(s)$.
There are many maps
$\phi_{st} \in {\rm PU}(2,1)$ such that
$\phi_{st}(E_s,p_s)=(E_t,p_t)$.
In this section, we construct a canonical one which
varies continuously with the parameters.

Let $\mu_s$ be the midpoint of the shorter of
the two arcs of $R$ which connects
$p_s$ to $q_s$.   Likewise define $\mu_t$.
Let $r_s=\Theta(E_S^*)=(e(s),\overline e(s))$.  Likewise define
$r_t$.
We define $\phi_{st}=\beta_{st} \circ \alpha_{st}$, where
$\alpha_{st}, \beta_{st} \in {\rm PU}(2,1)$ stabilize $R$ and
\begin{itemize}
\item[1.] $\alpha_{st}$ is a rotation such
that $\alpha(\mu_s)=\mu_t$,
\item[2.] $\beta_{st}$ is a pure translation along the line which
joins $\alpha_{st}(r_s)$ to $r_t$.
\end{itemize}
By construction, $\phi_{st}$
commutes with the map $r$, which stabilizes $R$.

Let $\rho$ be the round metric on $S^3$.   For any map
$\gamma \in {\rm PU}(2,1)$, define
$|\gamma|=\sup_{x \in S^3} \rho(x,\gamma(s))$.
In particular, define $\varepsilon_{st}=|\phi_{st}|$.
As $t \to s$, we have $\varepsilon_{st} \to 0$.

For any set $X \subset S^3$ let
$X^{\varepsilon}$ be the $\varepsilon$ tubulr
neighborhood of $X$, as measured in $\rho$.
Define
\begin{equation}
B(s,t)=(\partial \Delta)^{\varepsilon}; \hskip 20 pt
\varepsilon=\varepsilon_{st}.
\end{equation}
$B(s,t)$ is a thin solid torus which surrounds
$\partial \Delta$. \enddemo

\demo{{\rm 8.2.} Main construction}
Define $\Omega(s,t)=\Omega_{\Delta}(s,t) \cup d\Omega_{\Delta}(s,t)$.
Here
\begin{equation}
\Omega'_{\Delta}(s,t)=\phi_{st}(\Omega_{\Delta}(s)); \hskip 15 pt
\Omega_{\Delta}(s,t)=\Omega'_{\Delta}(s,t) \cup B(s,t).
\end{equation}
Finally, define
\begin{equation}
\Z(s,t)=\Xi \cup \Omega(s,t).
\end{equation}
For later use, we define the elevation map
\begin{equation}
\Psi_{st}=\Psi_s \cup \phi_{st}^{-1}.
\end{equation}
Here $\Psi_s$ is the $\zeta$-elevation map
associated to the $(E_s,p_s)$.
By construction $\Psi_{st}$ is
the $\zeta'$-elevation map for
$(E_t,p_t)$.  Here $\zeta'=\phi_{st}(\zeta)$.

By construction, 
$$\Omega'_{\Delta}(s,t)=
\Omega(E_t,p_t;\phi_{st} \partial \Delta).$$
Thus, $\Psi_{st}$ maps the hybrid sectors
in $\Omega'_{\Delta}(s,t)$ to standard position. 
In particular, the hybrid sectors in
$\Omega'_{\Delta}(s,t)$ are as well
adapted to the element $g_t$ as are the
hybrid sectors in the discarded $\Omega_{\Delta}(t)$.

At the same time $\Omega_{\Delta}'(s,t)$ is a
${\rm PU}(2,1)$-image of $\Omega_{\Delta}(s)$, and
thus has essentially the same geometry.
The object $\Z(s,t)$ retains the dihedral
symmetry of $\Z(s)$.  Finally, the set
$B(s,t)$ seals up the gaps created by the
discontinuous nature of our construction.

We will not actually deal with  $\Z(s,t)$ directly.
Rather, we will estimate the mismatch between
$\Z(s,t)$ and $\Z(s)$.

Define $\Omega'_d(s,t)=\phi_{st}(\Omega_d(s)).$
The remainder of the chapter is devoted
to proving the following result:\enddemo

\proclaimtitle{the year lemma}
\proclaim{Lemma}
\label{finite}
Let $a=.00046$ and
$b=.00221${\rm .}
Let $s=s_k$ and $t \in I_k${\rm .}
\begin{itemize}
\ritem{1.} $\varepsilon_{st}<a${\rm .}
\ritem{2.} $g_t(\Xi \cup d\Omega_{\Delta}(s,t)) \subset
 (g_s(\Xi \cup d\Omega(s))^b${\rm .}
\ritem{3.} $g_t(d\Omega'_d(s,t)) \subset
(g_s(d\Omega_d(s)))^b${\rm .}
\end{itemize}

\endproclaim

To avoid redundancy, we will just prove Statements 1 and 2.
Statement 3 has the same proof as Statement 2. 
\medbreak
 8.3. {\it Computational ingredients}.

\proclaim{Lemma}
\label{expp}
Let $s_k$ be as above and let $t \in I_k${\rm .}
There is a constant $a_k$ such that
\begin{itemize}
\ritem{1.} $8||e(s_k)|-|e(t)|| +
 |\arg e(s_k)-\arg e(t)|<a_k<a${\rm .}
\ritem{2.} $4 a_k+
2|s_j-t|<b$.
\end{itemize}
\endproclaim

\demo{Proof}
Let $\underline s=s_0,s_1, \ldots ,s_{730}=\overline s$ be the set of
$731$ maximally and evenly spaced points in $J_0$.
A straightforward calculation, detailed in Section~12, shows that,
for successive values $s_j, s_{j+1}$ we have

\begin{itemize}
\item[1.] $8(||e(s_j)|-|e(s_{j+1})||+10^{-5}) +  
 |\arg e(s_j)-\arg e(s_{j+1})+10^{-5}|=a_j<a.$
\item[2.] $4 a_j + 
2|s_j-s_{j+1}|<b.$
\end{itemize}

Each partition interval has the form
$[s_{k-1},s_{k+1}]$, with midpoint $s_k$.
Here $k$ is odd.
The two items above, combined with the
near monotonicity lemma, give the result.
\enddemo

Let $\Pi_Z: \C^2 \to \C$ be projection onto the first factor.
Referring to Equation \ref{As} define
\begin{equation}
\Sigma_k=\{x \in S^3|\ \sup_{s \in J} |\Pi_Z(x)-P_0(s)| \geq 1/k\},
\hskip 15 pt
P_0(s)=\frac{-\overline A_1(s)}{A_2(s)}.\qquad
\label{basicset}
\end{equation}

The Lipshitz lemma and the variation lemma below explain the significance of
this set.

\proclaimtitle{part of Estimate 1}
\proclaim{Lemma}
\label{part}
 For all $s \in J_0${\rm ,} $(\Xi \cup d\Omega_{\Delta}(s))^b
\subset \Sigma_3${\rm .}
\endproclaim

{\it Proof}.
See Section~9, Estimate 1.
\medbreak

For later use, 
let $F_k \subset \Sigma_k \times \Sigma_k$ denote those
points $(x,y)$ which can be joined by a geodesic
segment$-$i.e. an arc of a great circle$-$which is contained
entirely in $\Sigma_k$.

\demo{{\rm 8.4.} Proof of Statement $1$}
Statement 1 of Lemma \ref{expp}, together with the
following lemma, finishes the proof of Statement 1 of the
year lemma. \enddemo

\proclaim{Lemma}
\label{perturb}
$\varepsilon_{st} \leq |\arg(e(s))-\arg(e(t))|+
8 ||e(s)|-|e(t)||.$
\endproclaim

\demo{Proof}
Recall that $\phi_{st}=\beta_{st} \circ \alpha_{st}$.
Here $\alpha_{st}$ is a rotation which maps
the midpoint $\mu_s$ to the midpoint
$\mu_t$.  
Now, $\mu_s$ is contained in the
line through the origin which also contains
the point $\Theta(E_s^*)=(e(s),\overline e(s))$.
The corresponding statement holds
for $\mu_t$.
From this it follows that
$$
|\alpha_{st}| \leq
|\arg(e(s))-\arg(e(t))|.
$$
It remains to show that
$|\beta_{st}|
\leq 8 ||e(s)|-|e(t)||.$
From the near monotonicity lemma, and from
explicit evaluations at the endpoints of $J$,
\begin{equation}
\sup_{s \in J} \frac{\sqrt 2}
{2 |e(s)|^2-1} <8.
\label{star2}
\end{equation}
Let $\hat R$ be the real slice bounded by the 
round $\R$-circle $R$.   Note that
$\hat R=\{(z,\overline z)\} \cap \C\H^2.$
We put the Euclidean metric on $\hat R$.
By construction, the points
$$
(0,0); \hskip 15 pt
X_1=\alpha_{st}(e(s),\overline e(s)); \hskip 15 pt
X_2=(e(t),\overline e(t))
$$
are all contained in a single line $L$.
The element $\beta_{st}$ is a translation
along $L$ which maps
$p_1=\alpha_{st}(e(s),\overline e(s))$ to
$p_2=(e(t),\overline e(t)).$
Note that
$$
\|p_1\|=\sqrt 2 |e(s)|; \hskip 10 pt
\|p_2\|=\sqrt 2 |e(t)|; \hskip 10 pt
\|p_1-p_2\|=\sqrt 2 ||e(s)|-|e(t)||.  $$
Combining this information with
Lemma \ref{flowbound} and Equation \ref{star2}  gives  our estimate on
$|\beta_{st}|$.
\enddemo

\proclaim{{C}orollary}
\label{close}
For any point $x \in \Xi \cup d\Omega_{\Delta}(s_k,t)${\rm ,}
there is some $y \in \Xi \cup d\Omega_{\Delta}(s_k)$
such that $\rho(x,y)<a_k${\rm .}
\endproclaim

\demo{Proof}
Let $s=s_k$.
This follows from the fact that $\phi_{st}$
moves every point of $S^3$ less than $a_k$.
Hence, $d\Omega'_{\Delta}(s,t)$ is contained in
the $a_k$-neighborhood of
$d\Omega_{\Delta}(s)$.
By construction, $B(s,t)$ is contained in
the $a$-tubular neighborhood
of $\Xi \cup d\Omega_{\Delta}(s)$.
\enddemo

\demo{{\rm 8.5.} Proof of Statement  $2$}

\proclaimtitle{the Lipschitz estimate}
\proclaim{Lemma}
\label{lip}
For all $(x,y) \in F_k${\rm ,} and all
$k \geq 0${\rm ,}
$$\sup_{(x,y) \in F_k}
  \frac{\rho(g_s(x),g_s(y))}{\rho(x,y)} \leq
\frac{2k}{27} \sqrt{27+31k^2}.$$
In particular{\rm ,} the quantity on the right-hand side is less than $4$ for $k=3${\rm .}
\endproclaim

\demo{Proof}
From the formula in Section~3.2 we have
$$g(z,w)=\bigg(\frac {w}{A_2 z+\overline{A_1}},
\frac{A_1 z-A_2}{A_2 z + \overline{A_1}}\bigg)=(g_1(z,w),g_2(z,w)).
$$
Note that
$|A_2(s)| \geq 3 \sqrt 3/2$ for $s \in J$.
Let $P_0=P_0(s)$ be as in Equation \ref{basicset}.
We compute
\begin{eqnarray*}
|\partial_z g_1(s)|= \bigg| \frac{A_2 w}{(A_2 z - \overline A_1)^2}\bigg|
& \leq &\frac{1}{|A_2||z-P_0|^2} \leq \frac{2k^2}{3 \sqrt 3};\\ \noalign{\vskip4pt}
|\partial_z g_2(s)|
=\bigg|\frac{A_2^2-|A_1|^2}{(A_2 z - \overline A_1)^2}\bigg|
&\leq &\frac{1}{|A_2|^2|z-P_0|^2} \leq \frac{4k^2}{27};\\ \noalign{\vskip4pt}
|\partial_w(g_1)|=
\bigg|\frac{1}{|A_2 z + \overline A_1|}\bigg|
&\leq &\frac{1}{|A_2||z-P_0|} \leq
\frac{2k}{3 \sqrt 3};\\ \noalign{\vskip4pt}
\partial_w g_2&=&0.
\end{eqnarray*}

The bound in this lemma follows from integrating the
differential bound
$$
\frac{\|dg(v)\|}{\|v\|} \leq
\sqrt{|\partial_z g_1|^2
+|\partial_w g_1|^2
+|\partial_z g_2|^2
+|\partial_w g_2|^2}
$$
along the arc of the great circle joining the two relevant points.
\enddemo

\proclaimtitle{the variation estimate}
\proclaim{Lemma}
\label{var}
For all $s,t \in J$ and $x \in \Sigma_k${\rm ,}
$$\rho(g_s(x),g_t(s)) \leq \frac{2k}{3}|s-t|.$$
In particular{\rm ,} the quantity on the right\/{\rm -}\/hand side
is less than $2|s-t|$ for $k=3${\rm .}
\endproclaim

{\it Proof}.
Recall that
$$g=\Theta \circ \tilde g \circ \Theta^{-1}; \hskip 15 pt
\Theta(z_1,z_2,z_3)=(z_1/z_3,z_2/z_3); \hskip 15 pt
\Theta^{-1}(z,w)=(z,w,1).$$
Let $X=(z,w) \in S^3$.
Let $\tilde X=\Theta^{-1}(X)=(z,w,1)$.

Inspecting the formula for
$\tilde g$,  we see that
the third coordinate of $\tilde g(\tilde X)$ has
norm at least
$$
|\overline A_1 + A_2 z| \geq
|A_2|/k \geq \frac{3\sqrt 3}{2k}.
$$
Here we have used the same  bound on $|A_2|$ as
in the previous section.
Let $N_0(\varepsilon) \subset N_0$ denote those
points $(z_1,z_2,z_3)$ such that $|z_3| \geq \varepsilon$.
The above calculation shows that
$\tilde g_s(\tilde X) \subset N_0(\frac{3 \sqrt 3}{2k})$.
We compute the differential
$$
d\Theta=\left[ \matrix{
1/z_3 & 0 & -z_1/z_3^2 \cr
0& 1/z_3 & -z_2/z_3^2}\right].
$$
Using the fact that $|z_1|^2+|z_2|^2=|z_3|^2$,
we see that the sum of the squares of the matrix
entries of $d\Theta$ is $3/|z_3|^2$.
It follows that
$$
\sup_{y \in N_0(\varepsilon)}
\frac{\|d\Theta(V)\|}{\|V\|} \leq \sqrt 3 \varepsilon^{-1}.
$$
Here $V$ is any vector tangent to $y$.

Therefore, given $X \in \Sigma_k$,
the path $g_s(X)$ has speed at most
$\frac{2k}{3}$ times the speed of
the path $\tilde g_s(\tilde X)$, at corresponding points.
An exercise in calculus shows that
$|A_1'(s)|<.5$ and $|A_2'(s)|<.5$ for all
$s \in J$.
From this it follows that
$\tilde g_s(\tilde X)$ has
speed at most $1$.  Hence,
the path
$g_s(X)$ has speed at most
$\frac{2k}{3}$.   This is a reformulation of Lemma 8.7.
\hfill\qed\medbreak

Let $x \in \Xi \cup d\Omega_{\Delta}(s,t))$ be arbitrary.
By Corollary \ref{close}, there
is some point $y \in \Xi \cup d\Omega_{\Delta}(s)$ such
that $\rho(x,y)<a_k$.
By Lemma \ref{part}, and the triangle inequality, the points 
$x$ and $y$ are both at least $b-a>a$ from $S^3-\Sigma_3$.
Hence, $(x,y) \in F_3$.
From Lemmas \ref{lip} and \ref{var}
$$\frac{\rho(g_s(x),g_s(y))}{\rho(x,y)} \leq
4; \hskip 25 pt
  g_s(x)-g_t(x)<2|s-t|.
$$
From the second statement of Lemma \ref{expp}, and the triangle inequality:
$$\rho(g_t(x),g_s(y)) \leq
  \rho(g_t(x),g_s(x))+
  \rho(g_s(x),g_s(y)) \leq
4 a_j + 2|s-t| < b.$$

\section{The discreteness proof: General case}

9.1. {\it Five disks}.
The set of disks $W_{ABC}(s)$ can be defined, for each
parameter, just as in Section~7.  Here we define a single
collection $W_{ABC}$, which contains $W_{ABC}(s)$,
for all $s \in J$.   This collection is the union
of the five disks defined below.

\nonumproclaim{Lemma/Definition}
 Let $W_A(s)$ be the disk in $W$ such that
$A_s \subset \partial W_A(s)${\rm .}
Then $W_A(s) \subset W_A(s')$ if $s>s'${\rm .}
In particular{\rm ,}  
$W_A=W_A(\underline s)$ is defined so that
$W_A(s) \subset W_A$ for all
$s \in J${\rm .}
\endproclaim

\demo{Proof}
Here $A_s$ is the hyperbolic geodesic containing the fixed points
of $h_{+,s}$ and $h_{-,s}$, which move
together symmetrically along $\partial W_2$ as
$s$ increases. \phantom{SayingIloveyou}
\enddemo

\nonumproclaim{Lemma/Definition} Let $W_B^{\pm}(s)$ be the disk tangent to
$W_2$ at the fixed point $p_{\pm,s}$ of $h_{\pm,s}${\rm ,}
and also
tangent to $\partial W${\rm .}
Let $\sigma=(\underline s+\overline s)/2${\rm .}
Let $W_B$ be the smallest disk{\rm ,} concentric
with $W_B^{\pm}(\sigma)$ which contains
both $W_B^{\pm}(\underline s)$ and
$W_B^{\pm}(\overline s)${\rm .}
Then $h_{\pm,s}(W_2) \subset W_B^{\pm}$
for all $s \in J${\rm .}
\endproclaim

\demo{Proof}
By construction,
$h_{\pm,s}(W_2) \subset W_B^{\pm}(s)$.
Since the disks $W_B^{\pm}(s)$ are all isometric
they are all Euclidean isometric
to each other, and move monotonically through
the annulus $W_1$, we see that
$W_B^{\pm}(s) \subset W_B^{\pm}$ for all $s \in J$. \phantom{morethanwords}
\enddemo

\nonumproclaim{Lemma/Definition} Let $W_C^{\pm}(s)=h_{\pm,s} h_{\mp,s}(W_{\pm}^*)${\rm .}
Then $W_C^{\pm}(s) \subset W_C^{\pm}(s')$ if
$s<s'${\rm .}   If $W_C^{\pm}=W_C^{\pm}(\overline s)$
then $W_C^{\pm}(s) \subset W_C^{\pm}$ for all
$s \in J${\rm .}
\endproclaim

\demo{Proof}
We consider the case of $W_C^+$, the other case following
from symmetry. Now,
$h_{+,s}h_{-,s}(\{0\} \times \R^2)$ is determined
by its endpoints
$$
h_+h_-(0,1)=\bigg(0,\frac{A_1-A_2}{A_2+\overline A_1}\bigg),
\hskip 15 pt
h_+h_-(0,-1)=\bigg(0,\frac{-A_1-A_2}{A_2-\overline A_1}\bigg).
$$
These are the functions from Equation \ref{As}.
An exercise in calculus show that the real parts
(and hence the arguments) of the second coordinates of the above points are
monotone for $s \in J$.
Explicit evaluations at the endpoints
of $J$ show that
$0<\arg h_+h_-(0,1)<\arg h_+h_-(0,-1)<\pi/2,$
that $\arg h_+h_-(0,1)$ decreases with the parameter
and that $\arg h_+h_-(0,-1)$ increases.
Hence $W_C^+(s) \subset W_C^+(\overline s)$
for all $s \in J$.
\enddemo

\demo{{\rm 9.2.} All words but two}
Let $a$ and $b$ be as in the year lemma.

\proclaimtitle{Estimate 1}
\proclaim{Lemma}
For all $s \in J_0,$
\begin{itemize}
\ritem{1.} 
$\Pi_W((\Omega_{\Delta}(s))^a)  \cap W_{ABC}=\emptyset${\rm ,}
\ritem{2.} $(\Xi \cup d\Omega_{\Delta}(s))^b
\subset \Sigma_3${\rm ,} the set from Section~{\rm 8.}
\end{itemize}
\endproclaim

{\it Proof}.
See Section~12.
\hfill\qed

\proclaim{{C}orollary}
\label{proms}
$\Pi_W(\Omega_{\Delta}(s,t)) \cap W_{ABC}=\emptyset${\rm .}
\endproclaim

{\it Proof}.
The same argument as in Corollary \ref{close} implies that
$\Omega_{\Delta}(s,t) \subset (\Omega_{\Delta}(s))^a$.
The rest of the argument is immediate.
\hfill\qed\medbreak

The same symmetry as in
Section~7 imples that 
$\Pi_W(\Omega(s,t)) \cap W_{ABC}=\emptyset.$
As in Section~7, let $U_{1,t}$ be the component of
$S^3-\heartsuit(s,t)$ which contains
both fixed points of $h_{+,t}h_{-,t}$.
Let $U_{2,t}=d(U_{1,t})$.
Once we establish the following lemma,
the arguments in Section~7 apply {\it verbatim\/} to prove that
$h(U_{2,t}) \subset U_{1,t}$ for all $h \in H_t-H_{2,t}$.
The reader should bear in mind, while reading our proof,
that the thin set $B(s,t)$ plays an important
topological role in it.

\proclaim{Lemma}
$U_{1,t} \not = U_{2,t}${\rm .}
\endproclaim

{\it Proof}.
Assume this result is false.
Let $x$ be as Lemma \ref{distinct}.
The same argument in Lemma \ref{distinct},
combined with Lemma \ref{proms},
implies that there is a loop $\gamma \in S^3-\Omega(s,t)$
which links either $\partial \Delta$ or
$d\partial \Delta$, or both.  Assume, without
loss of generality, that it is the former.
$B(s,t)$ is a thin solid torus containing
$\partial \Delta$ as its core curve.
By construction,
$\eta=\partial \Omega'_{\Delta}(s,t) \subset B(s,t)$.
This latter curve runs nearly parallel to 
$\partial \Delta$, and is much longer than
the cross sectional diameter of $B(s,t)$.
Hence, $\eta$ is nontrivial in the first 
homology of $B(s,t)$.   In particular,
$\gamma$ and $\eta$ are linked.
The same argument as in Lemma \ref{distinct}
implies that $\gamma$ intersects
$\Omega'_{\Delta}(s,t)$, which is a contradiction.
\hfill\qed

\advance\eqcount by 57
\demo{{\rm 9.3.} Words of length $2$}
We found the following barrier function
$f_s:\break S^1 \to~\R$ experimentally:
\begin{equation}
f_s(u)=t(s)+(.21+.11(s-\overline s))\ \Im(\overline z(s) u); \hskip 15 pt 
(z^2(s),2t(s))=\Psi_s \circ g_s(\zeta).\enspace
\label{barrierfunction}
\end{equation}
Here $\zeta$ is as in Section~7.
We choose the branch of $\sqrt{}$ so that
$z(s)$ varies continuously.

Let $\Lambda(s)$ and $\Lambda_j(s)$
be the sets associated to the barrier
$(f_s,\Psi_s)$, as in Section~7.
Let $I_1(s)=(-\infty,t(s)]$ and
$I_2(s)=[t(s),\infty)$.\enddemo

\proclaimtitle{Estimate 2}
\proclaim{Lemma}
Let $c=a+b${\rm .}
For all $s \in J_0$ and for all
$t \in I_s${\rm ,}
\begin{itemize}
\ritem{1.} $\Psi_s(\Xi^c \cup (d\Omega_{\Delta}(s)^c)
\subset \Lambda_1(s).$
\ritem{2.} $\pi_{\R} \circ \Psi_s((d\Omega_d(s))^c) \subset 
I_1(s)${\rm .}
\ritem{3.} $\Psi_s(G_s)
\subset \Lambda_2(s)${\rm .}  Here
$G_s=(g_s(\Xi))^c \cup  (g_s(d\Omega_{\Delta}(s)))^c${\rm .}
\ritem{4.} $\pi_{\R} \circ \Psi_s(H_s) \subset I_2(s)${\rm .}
Here $H_s=(g_s(d\Omega_d(s)))^c${\rm .}
\ritem{5.} $\pi_{\R} \circ \Psi_{st}(g_t(p_t))>\max f_s${\rm .}
\ritem{6.} $g_t(p_t) \not \in (\Xi \cup d\Omega_{\Delta}(s))^c${\rm .}
\end{itemize}

\endproclaim

\demo{Proof}
See Section~12.
\enddemo

\proclaim{{C}orollary}
\label{main2}
Let $s=s_k \in J_0${\rm .}
For all $t \in I_k${\rm ,}  
\begin{itemize}
\ritem{1.} $\Psi_{st}(\Xi \cup d\Omega_{\Delta}(s,t))
\subset \Lambda_1(s).$
\ritem{2.} $\pi_{\R} \circ \Psi_{st}(d\Omega_d(s,t)) \subset 
I_1(s)${\rm .}
\ritem{3.} $\Psi_{st} \circ g_t(\Xi \cup d\Omega(s,t))
\subset \Lambda_2(s)${\rm .}
\ritem{4.} $\pi_{\R} \circ \Psi_{st} \circ g_t(d\Omega_d(s,t))
\subset I_2(s)${\rm .}
\ritem{5.} $\pi_{\R} \circ \Psi_{st}(g_t(p_t))>\max f_s${\rm .}
\ritem{6.} $g_t(p_t) \not \in \Xi \cup d\Omega_{\Delta}(s,t)${\rm .}
\end{itemize}

\endproclaim

\demo{Proof}
Statement 1 of the year lemma, combined with the
triangle inequality, implies that
$$\phi_{st}^{-1}(\Xi \cup d\Omega_{\Delta}(s,t))\subset
  (\Xi \cup d\Omega_{\Delta}(s))^{2a} \subset
  (\Xi \cup d\Omega_{\Delta}(s))^c.
$$
We also know that $\Psi_{st}=\Psi_s \circ \phi_{st}^{-1}$.
Hence, the set defined on the left-hand side of Statement 1 of the corollary
is contained in the set defined on the left-hand side of Statement 1 of
Estimate 2.
This proves Statement 1 of the corollary.
Statements 2--4 of the corollary have very similar proofs.
Statement 5 of the corollary has just been copied down
from Statement 5 of Estimate 2, for the sake of exposition.
Statement 6 of the corollary follows from Statement 6 of
Estimate 2, and Statement 2 of the year lemma.
\enddemo

Our goal is to show that
$\Z(s,t) \cap g_t(\Z(s,t))=\emptyset$ for
$t \in I_k$.   We   assume $t<\overline s$.
Once this is done, the same argument as in
Section~7 shows that
$(U_{1,t},U_{2,t})$ is a compressing
pair for $\rho_t$.

\proclaim{Lemma}
\label{star}
$\Z(s,t) \cap g_t(\Z(s,t) \subset \Omega_{\Delta}'(s,t) \cap
g_t(\Omega'_{\Delta}(s,t))${\rm .}
\endproclaim

\demo{Proof}
By Lemma \ref{nat},
$$\Psi_{st}(\Omega'_{\Delta}(s,t)-\{p_t\}) \subset
\Psi_{st}(\partial \Omega'_{\Delta}(s,t)) =
\Psi_s(\partial \Delta) \subset \Lambda_1(s).$$
Combining this with Statement 1 of 
Corollary \ref{main2}, we have
$$\Psi_{st}(\Z(s,t)-\{p_t\} ) \subset \Lambda_1(s).$$
Combining this with
Statements 2 and 6 of Corollary \ref{main2} we have
$$\Z(s,t) \cap g_t(\Xi \cup d\Omega_{\Delta}(s,t))=\emptyset.$$
The element $r$
conjugates $g_t$
to $g_t^{-1}$ and preserves every relevant set.
By symmetry, then,
$$g_t(\Z(s,t)) \cap (\Xi \cup d\Omega_{\Delta}(s,t))=\emptyset.$$
The lemma now follows from the definition
of $\Z(s,t)$.
\enddemo

\proclaim{Lemma}
$g_t(\Omega'_d(s,t)) \cap \Omega'_d(s,t)=\emptyset${\rm .}
\endproclaim

\demo{Proof}
First of all, we observe that
$g_t(p_t) \not = p_t$, and 
$\Omega'_d(s,t)  \cap E_t=\{p_t\}$. Also,
$g_t(\Omega'_d(s,t) ) \cap E_t=\{g_t(p_t)\}$.
Hence, 
$$p_t \not \in
g_t(\Omega'_d(s,t)) \cap \Omega'_d(s,t).
$$
By Lemma \ref{nat} and Statement 2 of Corollary \ref{main2},
$$
\pi_{\R} \circ \Psi_{st}(\Omega'_d(s,t)-\{p_t\})=
\pi_{\R} \circ \Psi_{st}(\partial \Omega'_d(s,t)) \subset
\pi_{\R} \circ \Psi_{st}(d\Omega_{\Delta}(s,t)) \subset I_1(s).
$$
Similarly, by Statement 4 of Corollary \ref{main2},
$$\pi_{\R} \circ \Psi_{st} \circ g_t(\partial \Omega'_d(s,t))
\subset \pi_{\R} \circ \Psi_{st} \circ g_t(d\Omega_{\Delta}(s,t))
\subset I_2(s).$$
Applying Statement 5 of Corollary \ref{main2} and
Lemma 3.4 to each of the foliating arcs
of $g_t(\Omega'_d(s,t)$, we get
$$
\pi_{\R} \circ \Psi_{st} \circ g_t(\Omega_d'(s,t)) \subset I_2(s).
$$
Since $I_1(s) \cap I_2(s)=\emptyset$,  the equations above
combine to prove this lemma.
\enddemo

\proclaim{Lemma}
\label{tame}
Let $s \in J_0$, the set defined in Section~{\rm 8.}
The hybrid disks
$\widetilde \Omega_h(s)$
and $\widetilde \Omega_v(s)$ are
tame{\rm .}
\endproclaim

\demo{Proof}
We check the
hypothesis of Lemma \ref{tamecrit} for
each parameter in $J_0$.
(Compare Equation \ref{tameevidence}.)
See Section~12 for details of this calculation.
\enddemo

\proclaim{Lemma}
\label{st6}
$\pi_{\R} \circ \Psi_s \circ g_s(\zeta)>\max f_s${\rm ,}
for all $s \in J_0${\rm .}
\endproclaim

\demo{Proof}
This is an explicit calculation.  See Section~12 for details.
\enddemo

\proclaim{Lemma}
$g_t(\Omega'_d(s,t)) \cap
\widetilde \Omega'_y(s,t)=\emptyset$ for
$y \in \{h,v\}${\rm .}
\endproclaim

\demo{Proof}
We will take $y=h$.  The case $y=v$ has the same proof.
We will show that the tame hybrid disks
$\widetilde \Omega_h(s,t)$ and the hybrid sectors
$g_t(\Omega'_d(s,t))$ satisfy the hypotheses of
Lemma \ref{disjo3}.

\begin{itemize}
\item[1.] The axis for $\widetilde \Omega'_h(s,t)$ is
$(E_t,p_t)$.   The axis for
$g_t(\Omega'_d(s,t))$ is
$(E_t,g_t(p_t))$.   Since $t \not = \overline s$,
and since the fixed points of $g_t$ are
not contained in the Clifford torus, we have
$q_t=g_t(p_t) \not = p_t$.    This shows
hypothesis 1.
\item[2.] Since $\partial \widetilde \Omega'_d(s,t) \subset
\Xi(s,t)$ and
$\partial g_t\Omega_d(s,t) \subset g_t(\Xi(s,t))$,
Lemma \ref{star} implies that
$\widetilde \Omega_d(s,t)$ and
$g_t(\Omega_h(s,t))$ satisfy hypotheses 2 and 3.
\item[3.] Note that $g_s(\zeta) \subset g_s(\partial \Omega_d(s))$.
Let $\alpha$ be the foliating arc of
$g_s(d\Omega_d(s))$ which contains $\zeta$.
Lemma \ref{st6} implies that
$\pi_{\R} \circ \Psi \circ g_s(\zeta)>\max f_s$.
Statement 5 of Estimate 2 implies that
$\pi_{\R} \circ \Psi \circ g_s(p_s)>\max f_s$.
Lemma \ref{monomap} now implies that
$\Psi_s(\alpha) \subset \Lambda_2(s)$.
From Lemma \ref{nat} and Estimate 2 we have
$\Psi_s(\Omega(E_s,p_s;C_s)-\{p_s\}) \subset \Lambda_1(s)$.
Since $p_s \not \in \alpha$, we see
that $\alpha \cap \Omega(E_s,p_s;C_s)=\emptyset$,
which is hypothesis 4.
\end{itemize}
\phantom{sun!}\vglue-50pt
\phantom{more sun}
\enddemo
 
\proclaim{Lemma}
$g_t(\widetilde \Omega_x(s,t)) \cap
\widetilde \Omega_y(s,t)=\emptyset$
for $x,y \in \{v,h\}${\rm .}
\endproclaim

{\it Proof}.
The proof is the same as in the previous case.
All we have to do is show that there is some
foliating arc of $g_s(\widetilde \Omega_y(s,t))$
which is disjoint from
$\widetilde \Omega_x(s,t)$.  This follows
from Estimate 2 and Corollary \ref{disjo2}.
\hfill\qed

\section{Numerical analysis}
\vglue-6pt
10.1 {\it The Clifford torus case}.
Recall that $\Xi$ is the part of the
dented torus which is also part of the Clifford torus; $\Xi$ is the ``outer rim" in Figure~10.1.
We subdivide $\Xi$ into five Clifford rectangles,
as shown in Figure 10.1.
(What looks like twelve rectangles is really $5$.)
Let $\{\Xi_j:\  j=1,2,3,4,5\}$ be the set of these pieces.
Let $\Xi(t,u)$ be the
rational parametrization of $\Xi_j$, as
described in Section~2.7.
\vglue-12pt
\figin{111.eps}{1000}
\vglue-9pt
\centerline{Figure 10.1}
\medbreak

Given a rectangle $Q \subset [0,1]^2$, let
$\Xi_j(Q)$ be the  Clifford subrectangle
of $\Xi_j$ parametrized by $Q$.
We now construct a ball
$B(\Xi_j,Q) \subset S^3$ such that
$\Xi_j(Q) \subset B(\Xi_j,Q)$.

Let $m$ be the midpoint of $Q$.
Let $V$ be the set of vertices of $Q$.
Let $B(\Xi_j,Q)$ be the
smallest ball centered at $\Xi_j(m)$ and  containing
the points $\Xi_j(V)$ in its boundary.
All four points are contained in the boundary
by symmetry.  Since each $\C$-arc foliating
$\Xi_j$ is contained in a semicircle,
we have $\Xi_j(Q) \subset B(\Xi_j,Q)$.

The {\it rectangle algorithm\/}
computes $B(\Xi_j,Q)$, returns
a $1$ if the
$t$-curves of $\Xi_j(Q)$ are longer than
the $u$-curves and returns a $2$ otherwise.
The return of $1$ tells the computer that
it can more efficiently decrease the
size of $\Xi_j(Q)$ by cutting $Q$ in half
so as to shorten the $t$-curves.  A
return of $2$ tells the computer the
opposite.

\demo{{R}emark}
The analysis above uses, in a crucial
way, the fact that the circular arcs
in question are contained in semicircles.
This property plays an important role in
our analysis
on several occasions.
The key fact about such arcs is that
their overall size is comparable to
the distance between their endpoints.
This is not true for circular arcs
in general.\enddemo

\advance\eqcount by 58

 10.2. {\it The hybrid sector case}.
Given an arc $\gamma \subset S^3$, let $E(\gamma)$ denote the
infimal $\varepsilon>0$ such that every point of
$\gamma$ is within $\varepsilon$ of either endpoint of $\gamma$.
For instance, if $\gamma=[0,1] \subset \R$ then
$E(\gamma)=1$.

Let $\Omega=\Omega(E,p;S)$ be a hybrid sector, with parametrization
$\Omega(t,u)$.   Let $T=[t_1,t_2]$ and $U=[u_1,u_2]$ be
subintervals of $[0,1]$.
Define
\begin{equation}
\Omega(T,u)=\bigcup_{t \in T} \Omega(t,u),
\hskip 25 pt
  \Omega(t,U)=\bigcup_{u \in U} \Omega(t,u),
\end{equation}
where $\Omega(t,U)$ is an $\R$-arc.
We call $\Omega(T,u)$ a {\it transverse arc\/}.
The goal of this chapter is to produce 
computable bounds $L_{\Omega}(t,U)$ and
$W_{\Omega}(T,U)$ such that
\begin{itemize}
\item[1.] $E(\Omega(t,U)) \leq
L_{\Omega}(t,U)$.
\item[2.] $E(\Omega(T,u)) \leq
W_{\Omega}(T,U)$  for any $u \in U$.
\end{itemize}
Our estimates are not general ones.
We have only verified them for the {\it clean subsectors\/}
of the hybrid sectors which actually
participate in Estimates 1 and 2.
(We will define below what we mean by {\it clean\/}.)

Let $T_1$ and $T_2$ be the two intervals obtained
from splitting $T$ in half.
Likewise define $U_1$ and $U_2$.
Let $m=(t_3,u_3)$ be the midpoint of $Q$.
We have
$$\rho(\Omega(t,u),\Omega(m) \leq
\rho(\Omega(t,u),\Omega(t_3,u))+
\rho(\Omega(t_3,u),\Omega(t_3,u_3)
$$ $$\leq\max (W(T_1,U),W(T_2,U))+
\max (W(t_3,U_1),W(t_3,U_2))=r(T,U).$$
Thus, $\Omega(Q)$ is contained in the ball
$B(\Omega,Q)$ centered at $\Omega(m)$ and
having radius
$r(T,U)$.

The {\it rectangle algorithm\/}
produces the ball
$B(\Omega,Q)$.  Also, it
returns a $1$ if
$L_{\Omega}(t_1,U) \leq W_{\Omega}(T,U)$, and
returns a $2$ if
$L_{\Omega}(t_1,U) \geq W_{\Omega}(T,U)$.
The integer values returned have a similar
purpose as in the previous section. 

\demo{{\rm 10.3.} Distortion of stereographic projection}
Define
\begin{equation}
U(\varepsilon)=\{(z,w) \in S^3|\ |z+1| \leq \varepsilon\}.
\end{equation}
Note that $U(2)=S^3$.
Let $\B_0$ be the map from Equation \ref{cayley}.
For any smooth map $T$, let $dT$ be the linear differential of
$T$.\enddemo

\proclaim{Lemma}
\label{disss}
If $V$ is a unit tangent vector to $S^3${\rm ,} based
at a point $p \in U(\varepsilon)${\rm ,} then
$\|d\B_0(V)\| \geq 1/(\varepsilon \sqrt 3)${\rm . } 
If $V$ is $\cal E$\/{\rm -}\/integral then
$\|d(\pi_{\C} \circ \B_0)(V)|\ \geq 1/\varepsilon${\rm .}
\endproclaim

\demo{Proof}
By rotational symmetry, continuity, and integration,
it suffices to assume that
$\B_0(p) \in \R \times \R$.
Note that $\B_0^{-1}(\R \times \R)$ consists of those
points $(z,w)$ such that $w/(z+1) \in \R$.

For $(z,w) \in \B_0^{-1}(\R^* \times \R)$, define
$$
\nu_r=\frac{(1+z)^2}{1+\overline z}
\left[\matrix{-\overline w \cr \overline z}\right]
; \hskip 15 pt
\nu_{\theta}=  \left[\matrix{0 \cr i (z+1)}\right]  ; \hskip 15 pt
\nu_t=  \frac{1}{2}\left[\matrix{i(z+1)^2 \cr iw(z+1)}\right].
\label{bootstrap}
$$
A vector $\nu$ is tangent 
(resp.\ $\cal E$-integral) to
$S^3$ at $p$ if and only if
$\Re\langle p,\nu\rangle=0$ (resp $\langle p,\nu \rangle=0$.)
Using this criterion, we see that all three vectors above are
tangent to $S^3$, and that
$\nu_r$ is also $\cal E$-integral.

Observe that all vectors have norm at most $\varepsilon$.
In the third case, this is a short exercise in calculus,
which uses the fact that $|z|^2+|w|^2=1$.
We will show that $d\B_0$ maps the vectors above
to an orthonormal basis.  This implies our
first differential bound.
We will see, in particular, that
$d(\pi_{\C} \circ \B_0)$ maps
$\nu_r$ to a unit vector.  This implies
the second bound.

Let $\beta_1$ and $\beta_2$ be the component functions of
$\B_0$.
We compute
\begin{eqnarray*}
d\beta_1&=&
\bigg(\frac{-w}{(1+z)^2},\frac{1}{1+z}\bigg);
\\ d\beta_1(\nu_r)=1; \hskip 15 pt
d\beta_1(\nu_{\theta})&=&i; \hskip 15 pt
d\beta_1(\nu_t)=0,
\\
d\beta_2(\nu_t)&=&
\pm \lim_{\varepsilon \to 0} \frac{1}{\varepsilon} \Im\bigg[ \frac{1-z}{1+z}-
\frac{1-z-i\varepsilon(z+1)^2/2}{1+z+i\varepsilon(z+1)^2/2}\bigg]=1.
\end{eqnarray*}
As we said above, $\nu_{\theta}$ is a real multiple of
$(0,iw)$, and this latter vector generates the flow
$(z,w) \to (z,\exp(it)\ w)$.   The map
$\B_0$ conjugates this flow to the flow on
$\cal H$ generated by $r\frac{d}{d\theta}$.
Thus, we see that $d\beta_2(\nu_{\theta})=0$.
The $\cal E$-integrality of $d\beta_1(\nu_r)$, together with
the fact that $d\beta_1(\nu_r)=1$, implies
that $d\beta_2(\nu_r)=0$.
\enddemo

\proclaim{Lemma}
\label{lippp}
Suppose $(z,w) \in S^3-U(\varepsilon)${\rm .}
If $V$ is a unit tangent vector based at $(z,w)$ then
$\|d\B_0(V)\| \leq \sqrt{\varepsilon^{-2}+5 \varepsilon^{-4}}
\leq \sqrt 6 \max(1,\varepsilon^{-2}).$
\endproclaim

{\it Proof}.
Let $\B_1$ be such that
$\B_0=\pi \circ  \B_1$. 
Here $\pi(z,w)=(z,\Im(w))$.
Note that
$\|d\B_0(V)\| \leq \|d\B_1(V)\|$.
We compute
$$
d\B_1=\left[\matrix{ 1/(1+z) & w/(1+z)^2 \cr
2/(1+z)^2 & 0}\right].
$$
The lemma follows straightaway.
\hfill\qed

\demo{{\rm 10.4.} Clean hybrid sectors}
Let $f_{\Omega}$ be the endpoint map for
$\Omega(E,p;S)$, as defined in Section~4.1.
We say that $\Omega(E,p;S)$ is
{\it clean\/} if
\begin{itemize}
\item[1.] $S$ is contained in a semicircle,
\item[2.] $f_{\Omega}(S) \subset E$ is contained in a semicircle,
\item[3.] $f_{\Omega}$ is injective on $S$.
\end{itemize}

Let $Y$ be a Heisenberg stereographic projection
mapping $\Omega$ to standard position.
$f_{\Omega}$ is injective on $S$ if and only if
$\pi_{\R} \circ Y$ is injective on $S$.
Now, $Y(S)$ is contained in
an ellipse, so that $\pi_{\R}$
has at most one local maximum and
one local minumum on $Y(S)$.

The {\it cleaning algorithm\/} does the following:
\begin{itemize}
\item[1.] Finds the points $\{p_i\}$ of $S$ where $f_{\Omega}$ is extremized.
\item[2.] Sorts the points $\{m\} \cup \{p_i\}$, according to the
order they appear on $S$.   Here $m$ is the midpoint of $S$.
\item[3.] Subdivides $S$ into subarcs $\{S_i\}$, bounded by the sorted points.
\item[4.] Checks that the hybrid sectors $\{\Omega(E,p;S_i\}$ are
clean.
\end{itemize}

We implement the first step using a simple iterative method.
We will explain how the method finds the point $\mu$ where
$\pi_{\R} \circ Y$ attains its minimum.  Finding the
maximum point is similar.
We first locate $\mu$ roughly, and thereby produce a 
map $\phi_0: [0,1] \to \widetilde S$ whose image contains
$\mu$.   Next, we evaluate the points
$\pi_{\R} \circ Y \circ (\alpha(j/3))$ for $j=0,1,2,3$.   The
minimum occurs either at $1/3$ or $2/3$.
In the first case, we define $\phi_1 = \phi_0 \circ A_1$.
Here $A_1$ is the affine map which takes $[0,1]$ to $[0,2/3]$.
In the other case, we define $\phi_1= \phi_0 \circ A_2$.
Here $A_2$ is the affine map which takes $[0,1]$ to $[1/3,1]$.
Iterating, we get a sequence of maps $\phi_0,\phi_1,\phi_2...$
whose images rapidly converge to $\mu$.

Let $\{S_i\}$ be the set of arcs produced by the
first three steps of the cleaning algorithm.
We will implement step 4 of the cleaning
algorithm by showing $p$
is farther from either of the two points of
$f_{\Omega}(\partial S_i)$ than these points
are from each other.

\medbreak 10.5. {\it Estimating the $\R$-arcs}.
Let $\Omega(E,p;S)$ be a hybrid sector.   Let $S_t$ be the
canonical parametrization for $S$. 
Let $A_t=f_{\Omega}(S_t)$ be the image of the
endpoint map, as in Section~5.3.
We define the {\it endpoint curves\/}:
\begin{equation}
E_p(t)=\B_0(-p \cdot A_t^{-1}); \hskip 25 pt
E_S(t)=\B_0(-S_t \cdot A_t^{-1}).
\end{equation}
We have used the group law on $S^3$ to define these
curves.
More generally, we define
\begin{equation}
E\Omega(t,u)=-\Omega(t,u) \cdot A_t^{-1}.
\end{equation}

\proclaim{Lemma}
\label{aff2}
$E\Omega(t,u)=u E_p(t) + (1-u) E_S(t).$
\endproclaim

{\it Proof}.
Apply Lemma \ref{aff} to the three points
$X=-p \cdot A_t^{-1}$ and
$Y=-S \cdot A_t^{-1}$ and
$Z=-A_t \cdot A_t^{-1}=(-1,0).$
\hfill\qed\medbreak

Combining Lemmas \ref{disss} and   \ref{aff2} we
see that, for all $u \in U$, we have 
$\rho(\Omega(t,u),\Omega(t,u_1) \leq L_1(t,U)$, where
\begin{equation}
L_1(t,U)=[2]\ |u_1-u_2|\ |\pi_{\C}(E_p(t))-\pi_{\C}(E_S(t))|.
\label{prelim}
\end{equation}
The $2$ above is bracketed for purposes which will
become clear in Section~10.7.

\demo{{\rm 10.6.} Estimating the transverse arcs}
A {\it dyadic interval\/} 
is an interval of the form
$[a/2^k,b/2^k]$, for $a,b,k$ integers.
When we run our main algorithm in the next section,
we will only have to deal with dyadic subintervals
of $[0,1]$.   Our analysis of transverse arcs
works for dyadic subintervals which have
length at most $1/4$.   If $T$ is a longer
interval, we set $W(T,U)=\pi$, the diameter of $S^3$.
Henceforth, we assume that $T$ is a 
dyadic subinterval of $[0,1]$, having length
at most $1/4$.
 \figin{112.eps}{1000}
 \centerline{Figure 10.6}
 \medbreak

We would like to compare 
$\omega_3=\Omega(t_3,u)$ and $\omega_4=\Omega(t_4,u)$, for
$t_3,t_4 \in T$.
The point $\omega_j$ lies in
the $\R$-arc $\Omega_j=\Omega(t_j,[0,1])$.
Let $A_j=A_{t_j}$ and $S_j=S_{t_j}$.
Let $R_j$ be right multiplication by $-A_j^{-1}$.
Let $X'_j=R_j(X_j)$ for any point $X_j \in \Omega_j$.
In particular, let $\Omega'_j=R_j(\Omega_j)$.
We have $A_1'=A_2'=(-1,0)$.
Since $R_j$ moves all points on $S^3$ the same distance,
and since
$A_1,A_2,A_3,A_4$ are contained, in order, on a semicircle, we get
\begin{equation}
\rho(\omega_3,\omega_4)<\rho(\omega'_3,\omega'_4)+\rho(A_3,A_4)
\leq \rho(\omega'_3,\omega'_4)+\rho(A_1,A_2).
\end{equation}

Let $\B_0$ be the stereographic
projection from Equation \ref{cayley}.
Note that $\B_0(0,-1)=\infty$.
Let $X''=\B_0(X')$ for all relevant $X$.
We equip $\cal H$ with the Euclidean norm.
In subsection 10.8 we will prove the following result.

\proclaim{Lemma}
\label{balls}
If $\Omega$ is a clean hybrid sector{\rm ,} produced when the
cleaning algorithm is applied to a hybrid sector
participating in Estimates $1$ and $2${\rm ,} then
$$\|p''_3-p''_4\| \leq .002+\|p''_1-p''_2\|, \hskip 30 pt
\|S''_3-S''_4\| \leq .002+ \|S''_1-S''_2\|.
$$
\endproclaim

By Lemma \ref{aff},
$\omega''_j=
u p''_j+ (1-u) S''_j$.
Combining this equation with
Lemma \ref{disss},
and Lemma \ref{balls}, we see that
$\rho(\omega_3',\omega_4') 
\leq W_1(T,u)$, 
where
\begin{eqnarray}
W_1(T,u)&=&
[2] \sqrt 3 \bigg(.002+ u \|\B_0(R_1(p_1))-\B_0(R_2(p_2))\|
 \\ 
&&+ \ (1-u) \|\B_0(R_1(S_1))-\B_0(R_2(S_2))\|\bigg). \nonumber
\end{eqnarray}
Again, $2$ above is deliberately bracketed.
The coefficients in the last equation are 
monotone for $u \in[0,1]$.  Hence, 
independent of the choice of $u \in U$, we get
$\rho(\omega'_3,\omega'_4) \leq W_1(T,U)$,
where
\begin{equation}
W_1(T,U)=\max(W_1(T,u_1),W_1(T,u_2)).
\end{equation}

\phantom{hoi}

\demo{{\rm 10.7.} Bootstrap argument}
Using the bounds above, we produce a ball $B'(\Omega,Q)$
such that $\Omega(Q) \subset B'(\Omega,Q)$.
Our calculation of the radius of $B'$  depends only on the
fact that $\Omega(Q) \subset U(2)$.
Given $B'$, we
find the best $\varepsilon>0$ such that $B' \subset U(\varepsilon)$.
Then, using
the full force of Lemma \ref{disss}, 
we can replace each occurrence of the bracketed
$2$ above by $\varepsilon$.
The smaller ball
$B(\Omega,Q)$ also contains $\Omega(Q)$.
This bootstrap argument reduces our computation time
roughly by a factor of $10$.
\enddemo

\demo{{\rm 10.8.} Proof of Lemma {\rm \ref{balls}}}
Let $\gamma: [0,1] \to \cal H$ be a curve.
Given a dyadic subinterval  $I \subset [0,1]$ let
$B(\gamma,I)$ be the ball in $\cal H$
which has, as a diameter, the segment joining the
two endpoints $\gamma(\partial I)$.
Let $I_m$ be the midpoint of $I$.
Say that $I$ is {\it good\/} if
$\gamma(I_m) \subset B(\gamma,I)$.

Let $\cal I$ denote the set of dyadic
subintervals of $T$.
Say that a {\it Lipschitz estimator\/} for
$\gamma$ is map $\Lambda: \cal I$$ \to \R$ such that
$\|\gamma(t)-\gamma(I_m)\|< \Lambda(I)$ for all
$t \in I$.   Consider this algorithm:

\begin{itemize}
\item[1.] Let LIST be the partition of $T$ consisting of just
$T$.
\item[2.] If LIST is empty, halt.
Otherwise...
\item[3.] Let $I$ be the first interval on LIST.
If $\Lambda(I)\leq \varepsilon$,
delete $I$ from LIST and go to
Step 2.   Otherwise...
If $I$ is good,
replace $I$ by the two intervals obtained by
cutting $I$ in half, and go to step 2.
Otherwise fail.
\end{itemize}

If this algorithm halts without failure, the information
constitutes a proof that $\gamma(t)$ is contained
within the $\varepsilon$ neighborhood of $B(\gamma,T)$.
The point is that there is a sequence
$T=T_0,T_1, \ldots ,T_k$ of dyadic intervals such
that $\|\gamma(t)-\gamma(\partial T_k)\|<\varepsilon$
and $B(\gamma,T_k) \subset B(\gamma,T_{k-1})...
\subset B(\gamma,T_0)$.

When we  run this algorithm, with the Lipschitz estimators
defined below, for $\varepsilon=.002$, and
for the curves $p''(t)=\B_0(-p \cdot A_t^{-1})$ and
$S''(t)=\B_0(-S_t \cdot A_t^{-1})$, it halts without
failure for every parameter in $J_0 \cup \{\overline s\}$.
See Section~12 for details of the calculation.
This establishes the lemma.
\enddemo

\demo{{R}emark}
For several parameters $s \in J_0$,
one of the hybrid sectors produced by the
cleaning algorithm is tiny.
In these cases, checking the ``control by balls''
condition involves the comparison of two tiny
quantities.   The roundoff error introduced by the
computer interferes with this comparison.
In these exceptional cases, we simply use our
Lipschitz estimators to verify, in advance,
that the relevant curves are entirely contained
in a ball of radius $.002$, thereby bypassing our
algorithm. \enddemo

\proclaim{Lemma}
Let $m$ be the midpoint of $[a,b]${\rm .}
Let
$$\lambda_1=\max(\rho(A_a,A_m),\rho(A_b,A_m)); \hskip 15 pt
\lambda_2=
|1+\pi_1(-p \cdot A_m^{-1})|
-\lambda_1 .$$
The quantity $\Lambda_p([a,b])=\lambda_1 \sqrt{6} \max(1,\lambda_2^{-2})$
is a Lipschitz estimator for $p''${\rm .}
\endproclaim

\demo{Proof}
Let $m$ be the midpoint of $[a,b]$.
Since $S|_T$ is contained in a semicircle, we have
$\rho(p_c,p_m) \leq \lambda_1$.
By the triangle inequality
$\rho(p'_c,p'_m) \leq \lambda_1$.
By definition, $p'_m=\B_0^{-1}(p''_m)$ and likewise
for $c$.  Therefore,
$$
\rho(\B_0^{-1}(p''(c)),\B_0^{-1}(p''(m)) \leq \lambda_1.
$$
By the triangle inequality,
$$
|\pi_1(\B_0^{-1}(p''(t)))+1|  \geq \lambda_2;
\hskip 15 pt t \in [c,m].
$$
The points $\B_0^{-1}(p''(a))$ and
$\B_0^{-1}(p''(c))$ can be connected by an arc
having length at most $\lambda_1$.   By Lemma \ref{lippp},
the restriction
of $\B_0$ to this arc expands distances along
this arc by at most $\sqrt 6 \max(1,\lambda_2^{-2})$.
\enddemo

Essentially the same argument gives

\proclaim{Lemma}
Let $m$ be the midpoint of $[a,b]${\rm .}
Let
\begin{eqnarray*}
\lambda_1&=&\max(\rho(A_a,A_m),\rho(A_b,A_m))+
\max(\rho(S_a,S_m),\rho(S_b,S_m)),\\
\lambda_2&=&
|1+\pi_1(-S \cdot A_m^{-1})|
-\lambda_1 .
\end{eqnarray*}
Then $\Lambda_p([a,b])=\lambda_1 \sqrt{6} \max(1,\lambda_2^{-2})$
is a Lipschitz estimator for $S''${\rm .}
\endproclaim

\section{Verifying the estimates}
\advance\eqcount by 66

All our estimates work, one parameter at a time.
Suppose the parameter is fixed.
The first step, in both cases,
is to split the relevant hybrid sectors into
clean hybrid sectors, and (for Estimate 2, Lemma 9.4) to split
$\Xi$ into $\{\Xi_j\}$, as in
Section~10.1.

Henceforth, we will say that a
{\it piece\/} is either a clean hybrid
sector participating in Estimates 1 (Lemma 9.1) and 2,
or one of the $\Xi_j$.
This chapter discusses how we verify that
Estimates 1 and 2 hold for an individual
piece.   We will also explain the 
remaining statements in Estimate 2, which do not pertain
to individual pieces.

\demo{{\rm 11.1.} The subdivision algorithm}
Let $Q_0=[0,1]^2$.
Given a {\it rectangle\/}
$Q=[a_1,b_1] \times [a_2,b_2] \subset  Q_0$, let
\begin{eqnarray}
\quad \Sigma_1(Q)&=&[a_1,(a_1+b_1)/2] \times [b_1,b_2] \cup
[(a_1+b_1)/2,a_2] \times [b_1,b_2],
\\
\Sigma_2(Q)&=&[a_1,a_2] \times [b_1,(b_1+b_2)/2] \cup
[a_1,a_2] \times [(b_1+b_2)/,b_2].
\end{eqnarray}

These are the two dyadic subdivisions of $Q$ into
two rectangles.
If $L=\{Q_k\}$ is a finite list of rectangles,
let $\Sigma_j(L)$ be the list obtained from $L$
by deleting the last rectangle $Q_n$ and appending 
$\Sigma_j(Q_n)$.

Let $\cal Q$ denote the set of rectangles of
$[0,1]^2$.
Say that a {\it subdivision test\/} is a map
$T: \cal Q$$ \to \{0,1,2\}$.
We are interested in proving the existence of a
partition of $Q_0$ 
into rectangles $\{Q_1, \ldots ,Q_n\}$ such that 
$T(Q_j)=0$ for $j=1, \ldots ,n$.   Here is the
{\it subdivision algorithm\/}: 
\begin{itemize}
\item[1.] Let LIST be a list with $Q_0$ as its only member.
\item[2.] If LIST is empty, halt.
If not, let $Q$ be the last rectangle in LIST.   
\item[3.] If $T(Q)=0$ delete $Q$ from LIST and go to step 2.
\item[4.] If $T(Q)=j \in \{1,2\}$ then replace LIST by $\Sigma_j$(LIST)
and go to step 2.
\end{itemize}
If the subdivision algorithm halts, the information
constitutes a proof that the desired partition
exists.\enddemo

\demo{{\rm 11.2.} Estimate $1$}
Let $a$ and $b$ be as in Estimate 1.
Let $\Omega$ be a piece.
Let $Q \subset [0,1]^2$ be a rectangle.
Let $P_0$ and $\Sigma_k$ be as in Equation
\ref{basicset}.
Here is the subdivision test: 

\begin{itemize}
\item[1.]
Let $B=B(\Omega,Q)$ be the ball
produced by the rectangle algorithm.
Let $(c,r)$ be the center and
radius of $B$.
Let $j \in \{1,2\}$ be the integer
produced by the rectangle algorithm.
\item[2.] For $i=1,2,3,4,5$, let
$\{(c_i,r_i)\}$ be the collections
of centers and radii for the disks in
$W_{ABC}$.  If
$\min_i (|\Pi_W(c)-c_i|-r-r_i) \leq a$,
let $T(Q)=j$.
Otherwise...
\item[3.] If $(\|\Pi_W \circ d(c))-P_0(s)\|-\frac{1}{3}-r\leq b$,
let $T(Q)=j$.  Otherwise let $T(Q)=0$.
\end{itemize}

$P_0(s)$ is evaluated from Equation \ref{basicset} (\S 8.3).
See Section~12.6 for a list of the centers and radii
of the disks used in this estimate.
\enddemo

\demo{{\rm 11.3.} Variation of the barrier function}
We use the notation of Section~9.3.
Let $f_s$ be the barrier function in 
Equation \ref{barrierfunction}.
Let $\pi_1: S^1 \times \R \to S^1$ be projection.
Let $\Psi_s$ be the elevation map associated to the parameter
$s$.  For $q \in S^3$ define
\begin{eqnarray}
F_s(q)&=&f_s \circ \pi_1 \circ \Psi_s(q)  + 
\pi_{\R} \circ \Psi_s(q),
\\
G_s(q)&=&\pi_{\R} \circ \Psi_s(q)-t(s).
\end{eqnarray}

We have defined this function so that
$$\begin{tabular}{cll}
{\rm 1.} & $\Psi_s(p) \in \Lambda_1(s)$ &\hbox{ if and only if 
$F_s(p) <0$.}\phantom{this is irritating but ok} \\
{\rm 2.} & $\Psi_s(p) \in \Lambda_2(s)$ &\hbox{ if and only if 
$F_s(p) >0$.} \\
{\rm 3.} & $\pi_{\R} \circ \Psi_s(p) \in I_1$ &\hbox{ if and only if  $G_s(p)<0$.} \\
{\rm 4.} & $\pi_{\R} \circ \Psi_s(p) \in I_2$ &\hbox{ if and only if 
$G_s(p)>0$.} 
\end{tabular}
$$

Let $B \subset S^3$ be a ball with center $c$ and radius $r$.
We define $F_s(B)=\bigcup_{x \in B} F_s(b)$.
If $p \in B$ we redefine $F_s(B)=[-10,10]$, which
causes $B$ to fail all further tests applied to it.
If $c=0$, we replace $B$ by the ball of radius
$r+.0000001$ centered at
$.0000001$.
Henceforth, we assume that $c \not = 0$ and $p \not \in B$.
We also assume that the parameter $s$ is fixed.
Let $\tilde c$ be the affinely normalized lift of $c$.

Given a vector $V=(z_1,z_2,z_3)$, let
\begin{equation}
\|V\|_2=\sqrt{|z_1|^2+|z_2|^2}.
\end{equation}
Define
\begin{equation}
L_+(B,V)=|\langle \tilde c,V\rangle| + r \|V\|_2; \hskip 15 pt
L_-(B,V)=\max(0,|\langle \tilde c,V\rangle| - r \|V\|_2).\quad
\end{equation}

Let $\hat p_s$ and $\hat q_s$ be any lifts of
the points $p_s, q_s \in E_s$.  Let
$\hat E_s^*$ be any vector polar to
$E_s$.  Let $r$ be the radius of $B$.  Define
\begin{eqnarray}
\alpha_1&=&\min \bigg(2 \pi,
\frac{r\|\hat E^*_s\|_2}{|L_-(B,\hat E^*_s)|} +
\frac{r\|\hat p_s\|_2}{|L_-(B,\hat p_s)|}\bigg),
\\
\varepsilon_1&=&\min(\alpha_1,2) \beta(s), \hskip 20 pt \beta(s)=.21+.11(s-\overline s)
, \\
\varepsilon_2&=& \frac{r\|\hat q_s\|_2}{L_-(B,p_s)}
+\frac{r\|\hat p_s\|_2 L_+(B,\hat q_s)}{L_-(B,p_s)^2} .
\end{eqnarray}

\proclaimtitle{the barrier test}
\proclaim{Lemma} 
With the above notation{\rm ,} $$\begin{tabular}{cll}
{\rm 1.} & \hbox{If $F_s(c)+\varepsilon_1+\varepsilon_2<0$} & \hbox{then
$\Psi_s(B) \subset \Lambda_1$.}\phantom{let go of our}\\
{\rm 2.} & \hbox{If $F_s(c)-\varepsilon_1-\varepsilon_2>0$} & \hbox{then
$\Psi(B) \subset \Lambda_2$.}\\
{\rm 3.} & \hbox{If $G_s(c)+\varepsilon_2<0$} & \hbox{then
$\pi_{\R} \circ \Psi_s(B) \subset I_1$.}\\
{\rm 4.} & \hbox{If $G_s(c)-\varepsilon_2>0$} & \hbox{then
$\pi_{\R} \circ \Psi_s(B) \subset I_2$.}\\
{\rm 5.} & \hbox{If $G_s(c)+\varepsilon_2+\beta(s)<0$} & \hbox{then
$\Psi_s(B) \subset \Lambda_1$.}\\
{\rm 6.} & \hbox{If $G_s(c)-\varepsilon_2-\beta(s)>0$} & \hbox{then
$\Psi_s(B) \subset \Lambda_2$.} 
\end{tabular}
$$
\endproclaim 

\demo{Proof}
Let $q \in B$ be an arbitrary point.
We define $(z_q,t_q)=\Psi_s(q)$.
Likewise $(z_c,t_c)=\Psi_s(c)$.
We first prove that
\begin{equation}
A=|\arg(z_c)-\arg(z_q)| \leq \alpha_1.
\label{argument}
\end{equation}
To see this, let $\gamma: [0,\delta] \to S^3$ be a
constant-speed geodesic which connects $c$ to $p$.
Let $\tilde \gamma$ be the affinely normalized lift of
$\gamma$.
We certainly have $A \leq 2 \pi$.
Define 
$X(t)=\arg( \Psi \circ \gamma(t))$.
By integration,
$$A \leq \sup_t \|X'(t)\| \delta <
\sup_t \|X'(t)\| r.$$
By definition,
$$X(t)=
\arg \bigg[\frac{\langle \tilde \gamma(t),\hat E_s^* \rangle}
     {\langle \tilde \gamma(t),\hat p_s \rangle}\bigg] =
\Im L(t); \hskip 15 pt 
L(t)=\log \bigg[\frac{\langle \tilde \gamma(t),\hat E_s^* \rangle}
     {\langle \tilde \gamma(t),\hat p_s \rangle}\bigg] .$$
Here $\log$ is a branch of the complex logarithm.
By calculus,
$$|X'(t)| \leq \|L'(t)\|=
\bigg| \frac{\langle \tilde \gamma'(t),\hat E^*_s \rangle}{\langle 
\tilde \gamma(t),\hat E^*_s \rangle}\bigg|+
\bigg| \frac{\langle \tilde \gamma'(t),\hat p_s \rangle}{\langle 
\tilde \gamma(t),\hat p_s \rangle}\bigg|.$$
The third coordinate of $\gamma'(t)$ is zero and
$\|\gamma'(t)\|=1$.
By the Cauchy-Schwarz inequality,
$$|\langle \gamma'(t),V\rangle| \leq \|V\|_2.$$
Here $V$ is any  vector whose third coordinate is $0$.
Another application of the Cauchy-Schwarz inequality
gives
$$|\langle \gamma(t),V \rangle| \geq L_-(B,V).$$ Therefore
$$|X'(t)| \leq \frac{\|\hat E^*\|_2}{L_-(B,\hat E^*)}+
\frac{\|\hat p\|_2}{L_-(B,\hat p)}= \varepsilon_1.$$
Equation \ref{argument} follows straight from this.

A very similar argument proves
\begin{equation}
|t_q-t_c| \leq \varepsilon_2.
\label{height}
\end{equation}
Given the definition of $f_s$, we have
\begin{equation}
\label{diffff}
\max f_s-\min f_s = 2\beta(s) ;\hskip 15 pt
|df_s/du| \leq \beta(s).
\label{derivatives}
\end{equation}
The first inequality in Equation 78 immediately
gives
\begin{equation}
\label{horbound1}
f_s \circ \pi_1 \circ \Psi_s(c)-f_s \circ \pi_1
\circ \Psi_s(q) \leq 2 \beta(s).
\end{equation}
Integrating the second inequality
in Equation 78 over the portion
of $S^1$ which contains $\pi_1 \circ \Psi_s(B)$,
we get
\begin{equation}
\label{horbound2}
f_s \circ \pi_1 \circ \Psi_s(c)-
f_s \circ \pi_1 \circ \Psi_s(q) \leq \alpha_1 \beta(s).
\end{equation}
Combining Equations \ref{height}, \ref{horbound1} and
\ref{horbound2} we see that
$|F_s(q)-F_s(c)|<\varepsilon_1+\varepsilon_2$.
Hence $F_s(q)<0$.
This is equivalent to the statement that
$\Psi(q) \subset \Lambda_1$.
The other statements have similar proofs.
\enddemo

\demo{{\rm 11.4.} Acting on the ball}
Let $\Sigma_k$ and $P_0(s)$, and $F_k$ be as in Lemma \ref{lip}.
Let $B \in S^3$ be a ball, with radius $r$ and center $c$.

Define
\begin{equation}
k=\frac{1}{|\pi_1(c)-P(s)|-r},
\hskip 20 pt K= \min(4,\frac{2k}{27} \sqrt{27+31k^2}).
\label{K}
\end{equation}
\enddemo

\proclaim{Lemma}
$g_s(B) \subset X(B,g_s)${\rm ,} where
$X(B,g_s)$ is the ball having center $g_s(c)$ and
radius $Kr${\rm .}
\endproclaim

\demo{Proof}
By construction,
every point $x \in B$ has the property that\break
$(c,x) \in F_k$.   Therefore, the Lipschitz
constant of $g|_B$ is at most $K$ by Lemma~\ref{lip}.   \phantom{sleet}
\enddemo

\demo{{\rm 11.5.} Confining a particular curve}
Let $s \in J_0$ and let
$I_s$ be the corresponding partition
interval.
Let $p_s$ be the point used in the definition of the
hybrid sectors of $\Omega_{\Delta}(s)$.\enddemo

\proclaim{Lemma}
\label{confine}
Let $\beta_s$ be the ball of radius $.001$ centered at
$g_s(p_s)${\rm .}  Then $(g_t(p_t))^a \in \beta_s$ for all 
$t \in I_s${\rm .}
\endproclaim

\demo{Proof}
Recall that $a=.00046$.
Let $P_0(s)$ be as in Equation \ref{basicset}.
Define
\begin{equation}
k=\frac{1}{|\pi_1(p_s)-P_0(s)|-2a-.02};
\hskip 15 pt K_1=\frac{2k}{27} \sqrt{27+31k^2};
\hskip 15 pt
K_2=\frac{2k}{3}|I_s|.
\end{equation}
From the year lemma, 
$\rho(p_s,p_t)<a$.  An exercise in calculus shows that
the real and imaginary parts of
$P_0(s)$ are monotone in $s$.   Evaluating at
the endpoints, we see that
$|P_0(s)-P_0(t)|<.02$  for all $s,t \in J$.  
By the triangle inequality,
$|\pi_1(p_t)-P_0(t)|<|\pi_1(p_s)-P_0(s)|+.02+a$.
We have chosen our constant $k$ so that the geodesic
in $S^3$ which joins $p_s$ to $p_t$ remains
in $\Sigma_k$.  Hence, $(p_s,p_t) \in F_t$.
Applying Lemma \ref{lip} we see that
$\rho(g_t(p_s),g_t(p_t))<K_1a$.     
Applying Lemma \ref{var} we see that
$\rho(g_s(p_s),g_t(p_s))<K_2$.     
By the triangle inequality,
$\rho(g_s(p_s),g_t(p_t)) \leq K_1 a + K_2$.
An explicit calculation, for all $s \in J$ shows
that the right-hand side of this equation is
at most $.00054=.001-a$.   See Section~12 for details.
\enddemo

{\elevensc {C}orollary 11.4.}
{\it Statement $6$ of Estimate $2$ {\rm (}\/from Lemma {\rm 9.4)} is true{\rm .}}
\medbreak

{\it Proof}.
Let $m_+(s)$ be the maximum value on $S_1$ attained
by $f_s$.   Let 
$$x_s=Y_s^{-1}(0,m_+(s)).$$
Here $Y_s$ is the Heisenberg stereographic
projection used to define $\Psi_s$.
Let $T_t=\pi_{\R} \circ \Psi_{st} \circ g_t$.
We compute explicitly that
$T_s(p_s)>\max f$ for all $s \in J_0$.
(See \S 12 for details of the calculation.)
We know that $g_t(p_t) \not = p_t$,
and $T_t(p_t)=\infty$.
Thus, $T_t(p_t) \not = \infty$.
To verify Statement 6, we just have to
rule out the possitility that
$T_t(p_t)=m_+(s)$.
If $T_t(p_t)=m_+(s)$ then
$
Y_{st} \circ g_t(p_t)=
Y_s \circ \phi_{st}^{-1} \circ g_t(p_t)=
Y_s(x_s)=(0,m_+(s)).
$
In other words, $\phi_{st}(x_s)=g_t(p_t)$.  
Hence, $x_s \in \beta_s$.
For each parameter $s \in J_0$ we compute explicitly
that $x_s \not \in \beta_s$.   See Section~12 for
details.
\hfill\qed

\demo{{\rm 11.6.} Estimate $2$}
Let $\Omega$ be a piece.
Let $Q \subset [0,1]$ be a rectangle. 
Here is the subdivision test for Estimate 2.

\begin{itemize}
\item[1.] Let $B=dB(\Omega,Q)$, where $B(\Omega,Q)$ is the ball returned from the
rectangle algorithm.  Let $j \in \{1,2\}$ be the
integer returned from the rectangle algorithm.
\item[2.] (Statement 1)
Use the Barrier Tests to see if $\Psi(B^a) \subset \Lambda_1$.
If not, set $T(Q)=j$.   Otherwise...
\item[3.]  (Statement 2, when relevant) Use the Barrier Tests
to see if $\pi_{\R} \circ \Psi(B^a)\break \subset I_1$.  
If not, set $T(Q)=j$.   Otherwise...
\item[4.] Let $B_g=X(g_s,B)$. 
\item[5.]  (Statement 3) Use the Barrier Tests to see if
$\Psi(B_g) \subset \Lambda_2$.  If not, let
$T(Q)=j$.  Otherwise...
\item[6.] (Statement 4, when relevant)  Use the Barrier Tests to
see if
$\pi_{\R} \circ \Psi(B_g)\break \subset I_2$.  If not, let
$T(Q)=j$.  Otherwise...
\item[7.] (Statement 5) If $B \cap \beta_s \not = \emptyset$ set
$T(Q)=j$.   Otherwise...
\item[8.]  Set $T(Q)=0$.
\end{itemize}
\enddemo

\vglue-24pt
\advance\eqcount by 82
\section{Implementation and roundoff error}

12.1. {\it Overview}.
The bulk of our computing is done in C.  We use
C code for
Lemma \ref{balls} and also for Estimates 1 and 2.
The C program 
only uses mathematical operations which are governed 
by the IEEE standards [I].   Namely,
\begin{itemize}
\item[1.] The $+,-,\times,\div$ and $\sqrt{}$ operations;
\item[2.] The $<$ and $>$ operations;
\item[3.] decimal $\Leftrightarrow$ binary conversion.
\end{itemize}
We implement our code using 
interval arithmetic, so that the computations
themselves produce the error bounds.

There are a number of places in the paper where
we need a small amount of fairly precise computation.
We perform these computations in Mathematica [W], which
has arbitrary precision arithmetic.
We compute all quantities to roughly
$20$ decimal places of precision.
While Mathematica is not guaranteed to be free
from computational bugs (and bugs have been
found!)
we think there is general agreement that
Mathematica does not have bugs in the basic operations
we use$-$the ordinary arithmetic operations
and the numerical extraction of cube roots.
We have ruled out such computational bugs by
random tests on our data$-$e.g. plugging the
root of a cubic back into the cubic.

The C program relies on some initial values,
depending on the parameter~$s$.   
For instance, the C program needs to have
the vector $E_s^*$ which in turn depends on
the eigenvalues of $g_s$.   Our Mathematica
code computes all initial values for the
C code, and stores them in auxiliary files.
When the C code runs, it reads in these files as needed.
Thus, the C code never has to compute any
quantity which depends on functions which
are not listed above.

We will first discuss certain features of our
Mathematica code.  Following this, we will
explain   how the computer represents
real numbers, and how we implement interval
arithmetic.  At the end of this chapter,
we give a record of all the calculations.
\advance\eqcount by 82

\demo{{\rm 12.2.} The Mathematica code}
Let $g_s$ be the matrix in Equation \ref{As}.
Let $E^*_s=(e(s),\overline e(s),1)$ be
the affinely normalized eigenvector associated to
$\lambda_s$.  
Equation \ref{theequation} gives $e(s)$ as a linear function
of $\lambda_s$.   Computing $e(s)$ accurately boils down to
computing $\lambda_s$ accurately.
Once we have an accurate value of $\lambda_s$, all
our Mathematica calculations are straightforward
implementations of what is explained in the paper.

Finding the roots of the charactistic
polynomial for $g_s$ leads to the formula
\begin{equation}
\lambda_s=\frac{a(s)+b(s)+c(s)}{12(1+s^2)},
\end{equation}
where
\begin{eqnarray}
a(s)&=&-68 + 64 i s - 4 s^2 ,
\\ \noalign{\vskip8pt} 
b(s)&=&8 \exp(4 \pi i/3) (-775 + 1734is +1203s^2-308is^3-69s^4
 \\  \noalign{\vskip6pt} 
&&+\ 6is^5 + s^6  
+6 \sqrt 3(1+s^2)^{3/2} \sqrt{3s^3-125})^{1/3},\nonumber\\ \noalign{\vskip8pt} 
 c(s)&=&\frac{5440-7936is -2688s^2 +256 is^3 +64s^4}{b(s)}.
\end{eqnarray}

This formula is valid, in particular, for $s \in J$.
Once Mathematica chooses branches
of the square and cube root functions, it can
compute the expression above
to, say, $1000$ decimal places.  We will stop
at $20$ decimal places.

Mathematica can consistently choose a branch of
the square root function for all $s \in J$.
Mathematica takes the positive branch of
$\sqrt{1+s^2}$ and the positive imaginary
branch of $\sqrt{3s^2-125}$.
Mathematica needs to take a cube root
in order to evaluate $b(s)$.
We have Mathematica
compute $b(s)$ along with $\lambda_s$.  
For successive computed parameters $s_1,s_2 \in J$, Mathematica checks
\begin{equation}
|b(s_j)|>50; \hskip 20pt
|b(s_1)-b(s_2)|<1.
\end{equation}
These bounds show that Mathematica cannot switch branches
of the cube root when computing successive parameters.
An explicit check shows that Mathematica chooses
the correct branches for $\overline s$.
Hence, the correct branches are always taken
during the computation.

Once again, we ruled out computational bugs in the
Mathematica code by randomly checking that the computed
eigenvalues and eigenvectors are indeed eigenvalues
and eigenvectors.
\enddemo

\demo{{\rm 12.3.} Performing the cleaning algorithm}
For each of the $366$ parameters $s \in J_0 \cup \{\overline s\}$,
we perform the cleaning algorithm on the hybrid
sectors in $\Omega_{\Delta}(s)$.
Each of the resulting clean hybrid sectors
$\Omega(E_s,p_s;S^j)$ is specified by
the parameter $s$, the
{\it axis},  which is the
collection $\{E^*_s,\hat p_s,\hat q_s,\hat \zeta\}$
normalized as in Section~5.4, and
an equalized choice of lifts
$\{\hat S_1^j,\hat S_2^j,\hat S_3^j\}$.
Here $S_1^j$ and $S_2^j$ are the endpoints of $S^j$ and
$S_3^j$ is the midpoint of $\widetilde S^j-S^j$. $\phantom{\sum^\int}$

The axis is the same for all the hybrid sectors.
The real and imaginary parts of the components of
all relevant vectors are what is stored in the file.
Everything is stored to $20$ significant digits.
(This is slightly more precision than is used in our C code.)
The routine {\tt Clean[data$s$,$s$]\/} performs the cleaning
algorithm for the parameter $s$ and stores the result
in the file {\tt data$s$\/}.   
The data contained in the created files is read in by the
C code at the time of its operation.
\enddemo

\demo{{\rm 12.4.} Doubles}
Our C code represents real numbers by
{\it doubles\/}.   
According to [I,  \S 3.2.2],
a double $[x]$ is an object of the form
$(s,e,f)$.  Here 
\begin{itemize}
\item[1.] $s$ is a single bit, determining the sign of $x$.
\item[2.] $e$ is an $11$ bit word, representing an integer
between $0$ and $2047$.
\item[3.] $f=\cdot b_1\cdots b_{52}$ is a  $52$    bit word.
\end{itemize}
The real number represented by $[x]$ is
\begin{equation}
r([x])=(-1)^s\ 2^{e-1023} \bigg(1+\sum_{j=1}^{52} 2^{-j} b_j\bigg).
\end{equation}
(The exceptions to this formula, which are
detailed in [I,  \S 3.2.2], do not arise
in our calculations.)
For example,
{\tt 3ff6a09e667f3bcd\/} represents the double closest to $\sqrt{2.0}$ 
The first bit of {\tt 3=0011\/} is {\tt 0\/}. 
Thus $s=0$.
The concatenation of {\tt 011\/} and
{\tt ff=11111111\/} is {\tt 01111111111\/}.
Thus $e=1023$, and
$(-1)^s 2^{e-1023}=1$.   The word {\tt 6a...\/} expands as
{\tt 01101010...\/}.   Thus
$$\sqrt{2.0}=1+  0+\frac{1}{4}+\frac{1}{8}+0+ \frac{1}{32}+0
+\frac{1}{128}+0\ldots =1.414\ldots \, .$$

To access the words comprising a double,
we introduce an {\it exact\/}:
\medbreak

\noindent
{\tt typedef struct \newline
$\{$int a,b;$\}$ exact.\/}\newline  
\newline
Given an exact {\tt E},  the
integer {\tt E.a\/} represents the first
$32$ bit word of the double (the {\it high bits\/}) and the 
integer {\tt E.b\/} represents the second
$32$ bit word of the double (the {\it low bits\/}).
Here is the code which enables one to convert
between doubles and exacts. 

\vfill
\noindent
{\tt exact double-to-exact(x)\/}
\newline
{\tt double x;\/}
\newline
$\{$ \newline
{\tt exact E;\/}  \newline 
{\tt double *pointer1;\/}  
\newline
{\tt int *pointer2;\/} \newline
{\tt pointer1=$\&$x;\/} \newline
{\tt pointer2=((int*) pointer1);\/} \newline
{\tt E.a=*pointer2;\/}  \newline
{\tt ++pointer2;\/}
\newline
{\tt E.b=*pointer2;\/}
\newline
{\tt return(E);\/}  \newline
$\}$
 \eject
\noindent 
{\tt double exact-to-double(E) \/} \newline
{\tt exact E;\/} \newline
$\{$ \newline
{\tt double x;\/} \newline
{\tt double *pointer1;\/} \newline
{\tt int *pointer2;\/} \newline
{\tt pointer1=$\&$x;\/} \newline
{\tt pointer2=((int*) pointer1);\/} \newline
{\tt *pointer2=E.a;\/} \newline
{\tt ++pointer2; \/} \newline
{\tt *pointer2=E.b; \/} \newline
{\tt return(x); \/} \newline
$\}$ 
\enddemo

\phantom{slow motion}

\demo{{\rm 12.5.} Interval arithmetic}
Let $\D$ be the set of doubles.
Define
$$\R_0=\{x \in \R|\ |x|<2^{1023}\}.$$
Summarizing [I,  \S 3.2.2, 4.1, 5.6], we see that there is a map
$\R_0 \to \D$, which maps each point $x \in \R_0$
to some $[x] \in \D$ which is closest to
$x$.  In case there are several equally close choices,
the computer chooses one, as detailed in\break [I,  \S 4.1].

Regarding the five basic operations,
[I,  \S 5] states that
{\it each of the operations shall be performed as if it
first produced an intermediate result correct to
infinite precision and with unbounded range, and
then coerced this intermediate result
to fit into the destination's format.\/}
Thus,
\vglue-4pt
\begin{equation}
\sqrt{[x]} := [\sqrt{r([x])}] \qquad
[x]*[y] := [r([x])*r([y])] ,\qquad  * \in \{+,-,\times,\div\}.  \quad
\end{equation}
\vglue12pt
We order exacts lexicographically, simply by concatenating
the two integers, treating them as a single
integer. 
The inclusion $r: \D \to \R$ induces a linear order 
on the doubles.
Beautifully, the linear order on positive doubles {\it coincides\/} with
the linear order we have imposed on the corresponding exacts.
Every double $x$ 
(except the largest and smallest) has a unique successor $x_+>x$ and
a unique predecessor $x_-<x$.   To compute these,
we convert $x$ to an exact, concatenate the integers,
add $1$ or $-1$, split apart the resulting integer, and
convert back.   

One exceptional case requires mention:
To avoid a certain kind of underflow error,
we redefined
$0_{\pm}$ to be the double closest to
$\pm 10^{-50}$.  In hindsight, this
redefinition seems unnecessary.  At
any rate, it is certainly harmless to our
calculations.

An {\it interval\/} is a pair $I=(x,y)$ of doubles,
such that $x \leq y$.
Say that $I$ {\it bounds\/} $z \in \R_0$ if
$x \leq [z] \leq y$.  This is true if and only if 
$x \leq z \leq y$.
Define
\vglue-6pt
\begin{equation}
[x,y]_o=[x_-,y_+].
\end{equation}
\vglue12pt\noindent 
Let $\partial_1 I$ and $\partial_2 I$
be the endpoints of an interval $I$.
For $j=1,2$, let $I_j=([x_j],[y_j])$ be
intervals.  For $* \in \{+,-,\times,\div\}$ we define
\vglue-6pt
\begin{equation}
I_1*I_2=(\min_{i,j}\partial_iI_1*\partial_jI_2,
\  \max_{i,j} \partial_iI_1*\partial_jI_2])_o;
\hskip 9 pt
\sqrt{I}=\left(\sqrt{[x]},\sqrt{[y]}\right)_o.\quad
\end{equation}
\vglue12pt\noindent 
 In all cases, one performs the relevant
computations on the endpoints of the interval
and then pushes the interval one click outward, to
guarantee that
\begin{itemize}
\item[1.] If $I_j$ bounds $x_j$ then $I_1*I_2$ bounds $x_1*x_2$.
\item[2.] If $I>0$ bounds $x>0$ then $\sqrt I$ bounds $\sqrt x$.
\end{itemize}
The only exception to this rule occurs when
we are computing $I_1/I_2$, and the
two endpoints of $I_2$ have different
signs. In this case we automatically
fail whatever computational test
we are working on.

A potential problem, related to the exceptional
case, is that we might
divide by a very small number, producing an
overflow error.  All the tests performed by
our code fail if a calculation produces a
sufficiently large number.  Thus, an overflow
error does not occur in the calculations
which are relevant to our proof.

We define a {\it complex interval\/} to be an expression
of the form $X+ i Y$, where $X$ and $Y$ are intervals.
We define a {\it vector interval\/} to be a
triple of complex intervals.   We define a 
{\it matrix interval\/} to be a triple of
vector intervals.  And so on. 
We say that a vector interval
bounds a vector if the components of the vector
interval bound the components of the vector.
Likewise for the other structures.

The algebra of the interval structures is identical
to the corresponding algebra for the usual structures.
At every step of our computation, the
actual object is bounded by the corresponding
interval version of the object.
Thus, if one of our algorithms halts with success,
as implemented,
the information constitutes a proof that a perfectly accurate
computing machine would also halt with success.
\enddemo

 12.6. {\it Record of the calculations}.
Our code is contained in the directory
{\tt~res/Computers/Proof\/}.
This directory is contained in the
University of Maryland College Park mathnet system.

The calculations from Section~3.4 are called
CALC1 and CALC2, and appear in the Mathematica code
in the same order they appear in Section~3.4.
The calculation CALC3 is the calculation
for the year lemma.
The calculation CALC4 establishes  
Lemma \ref{tame}.
The calculation CALC5 establishes
Lemma \ref{st6} and CALC6 establishes Lemma \ref{confine}.
Also, CALC7 establishes Corollary 11.4.
We performed the first three calculations on Sept.\ 2, 1999.
The remaining four were done on Sept.\ 30, 1999.
For all calculations, we used the
computer legendre.umd.edu, which is a Sparc Ultra~5.

Our Mathematica code also calculates
the five disks used in estimate 1.
We computed these disks on Sept.\ 2, 1999, on the same machine.
Here is the result:
\begin{itemize}
\item[1.] $W_A$ has center
$1.0757057484009542760...$ and\newline   radius
$.3964124835860459412...$.

\item[2.]  $W_B^{\pm}$ has 
center 
$.8425731633320585661... \pm
(.13646924571449026840...)i$ and
radius $.1513081547068024...$.

\item[3.] $W_C^{\pm}$ has center
$0.9316949906249123735...  \pm
(0.3892740451354294927...)i$ and radius
$.13996370162172745806...$
\end{itemize}

We performed the C calculations, during the period from
Sept.\ 24, 1999 to Sept.\ 28,  1999.  We used, in parallel,
the following Sparc Ultra 5 workstations: 

 \bigbreak
 {\tt
legendre.umd.edu 

noether.umd.edu 

kummer.umd.edu 

galois.umd.edu 

pascal.umd.edu

sylow.umd.edu 

descartes.umd.edu 

maclaurin.umd.edu 

monge.umd.edu 

poisson.umd.edu 

simpson.umd.edu 
\/}
\bigbreak
 \noindent
These ghosts of famous departed mathematicians
all took part in establishing our theorem.

%\vglue-6pt
\section{Computer plots}
%\vglue-6pt

Figure 13.1 shows $\Pi_W(d\Omega_d(\overline s))$.   The 
five grey disks are $W_{ABC}$.   The curves in the
picture are individual foliating $\R$-arcs, which attach to
$\partial W_2$.  

 Figure 13.2 shows a close up of Figure 13.1, with
 some additional points plotted.

 \phigin{1.eps}{250}
 \centerline{Figure 13.1}
 \vfill
\phigin{2.eps}{250}
\centerline{Figure 13.2}
 \eject
Figure 13.3 shows $\Pi_W(d\Omega_h(\overline s))$.   The 
five grey disks are $W_{ABC}$.   The curves in the
picture are individual foliating $\R$-arcs, which attach to
$\partial W_2$.
 \phigin{3.eps}{250}
\centerline{Figure 13.3}
 \vfill
Figure 13.4 shows a close up of Figure 13.3, with
some additional points plotted.

\phigin{4.eps}{250}
\centerline{Figure 13.4}
  
\eject
Figure 13.5 shows $\Pi_W(d\Omega_h(\overline s))$.   The 
five grey disks are $W_{ABC}$.   The curves in the
picture are individual foliating $\R$-arcs, which attach to
$\partial W_2$.

\phigin{5.eps}{250}
\centerline{Figure 13.5}
\medbreak

Figure 13.6 shows a close up of Figure 13.6, with
some additional points plotted.

\phigin{6.eps}{250}
\centerline{Figure 13.6}
\medbreak

Figures 13.7 shows 
$\Psi_{\overline s}(d\Omega_d(\overline s))$ and
$g_{\overline s}(\Psi_{\overline s}(d\Omega_d(\overline s)))$.
The 
cylinder $\Upsilon$ has been identified to
$\R/2 \pi \Z \times \Z$ via
$(u,t) \to (\arg u,t)$.
The sine curve is the common boundary
of $\Lambda_1(f)$ and $\Lambda_2(f)$, with
$f$ as in Equation \ref{specialf}.

\phigin{7.eps}{250}
\centerline{Figure 13.7}
\medbreak

Figure 13.8 shows 
$\Psi_{\overline s}(d\Omega_v(\overline s))$ and
$g_{\overline s}(\Psi_{\overline s}(d\Omega_v(\overline s)))$.
When $h$ replaces $v$, the picture looks exactly like
Figure 13.8, but
rotated $180^\circ$  about the points of symmetry of
$\Lambda(\overline s)$.

\phigin{8.eps}{250}
\centerline{Figure 13.8}
\medbreak

In Section~10 we  break $\Xi$ into $5$ pieces:
$\Xi=\Xi_1 \cup ... \cup \Xi_5$.
(Compare Figure 10.1.)
Figure 13.10 shows
$\Psi_{\overline s}(\Xi_2 \cup \Xi_3)$ and
$\Psi_{\overline s}(g_{\overline s}(\Xi_2 \cup \Xi_3) )$.
\vglue12pt
\phigin{11.eps}{250}
\centerline{Figure 13.10}
\bigbreak

Figure 13.11 shows a closeup of Figure 13.10.
A somewhat different set of points has been
plotted.
\phigin{12.eps}{250}
\centerline{Figure 13.11}
\medbreak
 
Figure 13.12 shows
$\Psi_{\overline s}(\Xi_1)$ and
$\Psi_{\overline s}(g_{\overline s}(\Xi_1))$.

\phigin{10.eps}{250}
\centerline{Figure 13.12}
\medbreak

We have $r(\Xi_2 \cup \Xi_3)=
\Xi_4 \cup \Xi_5$.   Thus, the picture for
$\Psi_{\overline s}(\Xi_4 \cup \Psi_5)$ and
$\Psi_{\overline s}(g_{\overline s}((\Xi_4 \cup \Psi_5))$ is
obtained by rotating Figures 13.10 and 13.11 about the point
of symmetry on $\Lambda_{\overline s}$

  \bye